\documentclass[12pt]{article}
\usepackage{amsmath,amssymb,amsthm}
\usepackage{comment}
\usepackage{color}
\numberwithin{equation}{section}
\newtheorem{thm}{Theorem}[section]
\newtheorem{prop}[thm]{Proposition}

\newtheorem{exam}[thm]{Example}
\newtheorem{rem}[thm]{Remark}
\newtheorem{lem}[thm]{Lemma}

\newtheorem{assum}[thm]{Assumption}
\newcommand{\pf}{{\em Proof.}}

\def\Capa{\mathop{\rm Cap}}

\makeatletter
\def\address#1#2{\begingroup
\noindent\parbox[t]{16cm}{%
\small{\scshape\ignorespaces#1}\par\vskip1ex
\noindent\small{\itshape E-mail address}%
\/: #2\par\vskip4ex}\hfill%
\endgroup}%

\title{Bottom crossing probability for symmetric jump processes}
\author{Yuichi Shiozawa
\footnote{Supported in part by the Grant-in-Aid for Scientific Research (C) 26400135.}}
\begin{document}
\maketitle 

\begin{abstract}
We determine the decay rate of the bottom crossing probability 
for symmetric jump processes under the condition on heat kernel estimates. 
Our results are applicable to symmetric stable-like processes 
and stable-subordinated diffusion processes on a class of 
(unbounded) fractals and fractal-like spaces. 
\end{abstract}

\section{Introduction}
We are concerned with the quantitative characterizations of 
{\it transience} and (non-point) {\it recurrence}  
for symmetric jump processes generated by regular Dirichlet forms. 
Such characterizations are expressed in terms of {\it lower rate functions}.
In this paper, we discuss the long time asymptotic estimates 
of the bottom crossing  probability 
related to lower rate functions.

For $\alpha\in (0,2]$, let 
${\bf M}=(\{X_t\}_{t\geq 0}, \{P_x\}_{x\in {\mathbb R}^d})$ 
be the symmetric $\alpha$-stable process on ${\mathbb R}^d$.
If ${\alpha=2}$, then ${\bf M}$ is the Brownian motion on ${\mathbb R}^d$. 
As is well known, if $d>\alpha$, 
then ${\bf M}$ is transient in the sense that 
the particle escapes to infinity eventually with probability one. 
On the other hand, if $d=\alpha \ (=1 \ \text{or} \ 2)$, 
then ${\bf M}$ is non-point recurrent in the sense that 
the particle comes arbitrarily close to the origin but never hits it 
with probability one. 
We can characterize these two properties quantitatively as follows: 
Assume that $d\geq \alpha$.  If $g(t)$ is a positive decreasing function on $(0,\infty)$ such that 
$g(t)\rightarrow 0$ as $t\rightarrow\infty$, then the function $r(t)=t^{1/\alpha}g(t)$ satisfies 
\begin{equation}\label{eq:0-1}
P\left(\text{there exists $T>0$ such that $|X_t|\geq r(t)$ for all $t\geq T$}\right)=1 \ \text{or} \ 0
\end{equation}
according as 
\begin{equation}\label{eq:test-stable}
\int_{\cdot}^{\infty}h_g(t)\frac{{\rm d}t}{t}<\infty \ \text{or} \ =\infty
\end{equation}
for $P=P_0$ and 
$$
h_g(t)=\begin{cases}
g(t)^{d-\alpha} & (d>\alpha), \\
\displaystyle \frac{1}{|\log g(t)|} & (d=\alpha).
\end{cases}
$$
When $\alpha=2$,  
Dvoretzky and Erd\"os \cite{DE51} and Spitzer \cite{Sp58} established 
this integral test for $d\geq 3$ and $d=2$, respectively (see also \cite[4.12]{IM74}).  
When $0<\alpha<2$, 
J. Takeuchi \cite{T64} and  J. Takeuchi and S. Watanabe \cite{TW64} obtained 
the test for $d>\alpha$ and $d=\alpha=1$, respectively.
If the probability in (\ref{eq:0-1}) is one, 
then the function $r(t)$ is called a \textit{lower rate function} of ${\bf M}$. 
This function expresses how fast the particle escapes to infinity for $d>\alpha$, 
and how arbitrarily close it comes to the origin for $d=\alpha$.
We can regard this function as the bottom of ${\bf M}$ for all sufficiently large time. 

Wichura \cite{W79} (see also \cite{W73} in the Brownian case) further 
proved that if the integral in (\ref{eq:test-stable}) is convergent,  
then there exists a positive constant $L_{d,\alpha}$ such that 
\begin{equation}\label{eq:w-type}
P\left(\text{$|X_t|<r(t)$ for some $t>T$}\right)
=(1+o(1))L_{d,\alpha}\int_T^{\infty}h_g(s)\frac{{\rm d}s}{s} \quad (T\rightarrow\infty)
\end{equation}
under some additional condition on the function $g(t)$. 
This equality gives the precise asymptotic behavior of the bottom crossing probability 
and related it to the integral in (\ref{eq:test-stable}).

The integral tests on lower rate functions are  extended to 
more general symmetric diffusion processes (see, e.g., \cite{BS05, I78, G99}) 
and symmetric jump processes (see, e.g., \cite{H70, K97, S16, SW15}). 
Among them, the full heat kernel estimates are utilized in \cite{BS05, I78, SW15} 
to establish zero-one law type results.
Our purpose in this paper is to determine 
the decay rate of the bottom crossing probability 
for a class of symmetric jump processes with no scaling property 
(see Theorems \ref{asymp-thm} and \ref{asymp-thm-critical} below).

Our approach here is based on that of Wichura \cite{W79}.  
However, the scaling property and the rotation invariance of 
symmetric stable processes on ${\mathbb R}^d$ played a crucial role 
in his proof.  
Instead of these properties, 
we make use of the full heat kernel estimates by following \cite{KKW15} and \cite{SW15}. 
Our results are applicable to 
symmetric stable-like processes (see \cite{C09,CK03,CK08})  
and a class of symmetric jump processes  
on (unbounded) fractals and fractal-like spaces (see Section \ref{sect:exam} below for details).

The rest of this paper is organized as follows:
In Section \ref{sec:result}, we first make assumptions on heat kernels 
and introduce the notion of lower rate functions. 
We then state our main results in this paper. 
In Section \ref{sect:exam}, we first calculate the decay rate of the bottom crossing probability 
for some lower rate functions. 
We then give examples to which our main results are applicable. 
In Section \ref{sect:hit}, we give estimates on the hitting time 
distributions of a process to closed balls during finite time interval.   
Using these estimates, we prove our result for the transient case in Section \ref{sect:proof}. 
The proof for the non-point recurrent case is given in Appendix \ref{sect:appendix-proof}
because this proof  is similar to that for the transient case.
Appendix \ref{sect:sub} is devoted to the calculation of Dirichlet forms and heat kernels 
for a class of subordinated diffusion processes, which will be mentioned in Section \ref{sect:exam}. 

Throughout this paper, the letters $c$ and $C$ (with subscript) 
denote finite positive constants which may vary from place to place. 
For positive functions $f(t)$ and $g(t)$ on $(1,\infty)$, we write 
$f(t)\asymp g(t) \ (t\rightarrow\infty)$ if there exist positive constants $T$, $c_1$ and $c_2$ such that 
$c_1g(t)\leq f(t)\leq c_2g(t)$ for all $t\geq T$.  
For nonnegative functions $f(x)$ and 
$g(x)$ on a space $S$, we write $f(x)\lesssim  g(x)$ (or $g(x)\gtrsim f(x)$) 
if there exists a positive constant $c$ such that 
$f(x)\leq cg(x)$ for all $x\in S$. 
We also write $f(x)\asymp g(x)$ if $f(x)\lesssim g(x)$ and $g(x)\lesssim f(x)$.

\section{Results}\label{sec:result}
We first recall the notion of Dirichlet forms from \cite{CF12} and \cite{FOT11}. 
Let $(M,d)$ be a locally compact separable metric space 
and $m$ a positive Radon measure on $M$ with full support. 
We write $C(M)$ for  the totality of continuous functions on $M$ 
and $C_0(M)$ for that of continuous functions on $M$ with compact support.  
Let $({\cal E}, {\cal F})$ be a Dirichlet form on $L^2(M;m)$, 
that is, $({\cal E}, {\cal F})$ is a closed Markovian symmetric form on $L^2(M;m)$. 
In what follows, we suppose that $({\cal E}, {\cal F})$ is {\it regular}:  
${\cal F}\cap C_0(M)$ is dense both 
in $L^2(M;m)$ with respect to the norm $\|u\|_{{\cal E}_1}=\sqrt{{\cal E}_1(u,u)}$, 
and in $C_0(M)$ with respect to the uniform norm $\|\cdot\|_{\infty}$. 
Here 
$${\cal E}_1(u,u)={\cal E}(u,u)+\|u\|_{L^2(M;m)}^2, \quad u\in {\cal F}.$$

Let ${\cal O}$ be the family of all open subsets of $M$. 
For $A\in {\cal O}$, we let 
$${\cal L}_A=\left\{u\in {\cal F} \mid \text{$u\geq 1$ \ $m$-a.e.\ on $A$}\right\}$$
and 
$$\Capa(A)=\begin{cases}
\inf_{u\in {\cal L}_A}{\cal E}_1(u,u), & \text{if ${\cal L}_A\ne \emptyset$},  \\
\infty, & \text{if ${\cal L}_A\ne\emptyset$}.
\end{cases}$$
For any  $A\subset M$, we define the {\it {\rm (1-)}capacity} of $A$ by  
$$\Capa(A):=\inf_{B\in {\cal O}, \, A\subset B}\Capa(B).$$
For  $A\subset M$, a statement depending on $x\in A$ is said to hold q.e.\ on $A$ 
if there exists a set $N\subset A$ of zero capacity such that the statement holds 
for every $x\in A\setminus N$. 
Here q.e.\ is an abbreviation for {\it quasi everywhere}.

We  suppose that the Beurling-Deny expression of $({\cal E}, {\cal F})$ 
(see \cite[Theorem 3.2.1, Lemma 4.5.4]{FOT11}) is given by 
\begin{equation}\label{b-d}
{\cal E}(u,v)=\iint_{M\times M\setminus{\rm diag}}(u(x)-u(y))(v(x)-v(y))\,J({\rm d}x{\rm d}y) 
\quad \text{for $u,v\in {\cal F}\cap C_0(M)$},
\end{equation}
where
${\rm diag}=\{(x,y)\in M\times M \mid x=y\}$ and 
$J({\rm d}x{\rm d}y)$ is a symmetric positive Radon measure on $M\times M\setminus{\rm diag}$.
We call $J$ the {\it jumping measure}  associated with $({\cal E}, {\cal F})$.

We write ${\cal B}(M)$ for the family of all Borel measurable subsets of $M$. 
Let $M_{\Delta}=M\cup\{\Delta\}$ be the one point compactification of $M$ 
and ${\cal B}(M_{\Delta})={\cal B}(M)\cup\{B\cup\{\Delta\} \mid  B\in {\cal B}(M)\}$. 
Let ${\bf M}=\left(\{X_t\}_{t\geq 0}, \{P_x\}_{x\in M}\right)$ be 
an $m$-symmetric Hunt process on $M$ generated by $({\cal E}, {\cal F})$ (\cite[Theorem 7.2.1]{FOT11}). 
A set $B\subset M_{\Delta}$ is called {\it nearly Borel measurable} 
if for any probability measure $\mu$ on $M_{\Delta}$, 
there exist $B_1, B_2\in {\cal B}(M_{\Delta})$ such that $B_1\subset B\subset B_2$ and 
$$P_{\mu}(\text{$X_t\in B_2\setminus B_1$ for some $t\geq 0$})=0.$$

By  \cite[Theorem 4.2.1]{FOT11}, a set $N\subset M$ is of zero capacity 
if and only if $N$ is {\it exceptional},  that is, 
there exists a nearly Borel measurable set $\tilde{N}\supset N$ such that $P_m(\sigma_{\tilde{N}}<\infty)=0$.  
Here $\sigma_{\tilde N}=\inf\{t>0 \mid X_t\in \tilde{N}\}$ is the hitting time of ${\bf M}$ to $\tilde{N}$.
We say that a set $N\subset M$ is {\it properly exceptional} 
if $N$ is nearly Borel measurable  such that $m(N)=0$ and $M\setminus N$ is {\it ${\bf M}$-invariant}, that is, 
$$P_x(\text{$X_t\in (M\setminus N)_{\Delta}$ and $X_{t-}\in (M\setminus N)_{\Delta}$ for any  $t>0$})=1 
\quad \text{for any $x\in M\setminus N$.}$$
Here $(M\setminus N)_{\Delta}=(M\setminus N)\cup \{\Delta\}$ and  $X_{t-}=\lim_{s\uparrow t}X_s$. 
Note that any properly exceptional set $N$ is exceptional and thus $\Capa(N)=0$ by \cite[Theorem 4.2.1]{FOT11}.

For $x\in M$ and $r>0$, let 
$B(x,r)=\left\{y\in M \mid d(y,x)<r\right\}$ 
be an open ball with radius $r$ centered at $x$ 
and let $V(x,r)=m(B(x,r))$. 
We make the following three assumptions:

\begin{assum}\label{assum:volume}
For any $x\in M$ and $r>0$, $V(x,r)<\infty$. 
Moreover, there exist constants 
$c_1\in (0,1]$, $c_2\in [1,\infty)$, $d_1>0$ and $d_2\geq d_1$ 
such that  
\begin{equation}\label{eq:volume-growth}
c_1\left(\frac{R}{r}\right)^{d_1}\leq \frac{V(x,R)}{V(x,r)}
\leq c_2\left(\frac{R}{r}\right)^{d_2} \quad 
\text{for all $x\in M$ and  $0<r<R$}.
\end{equation}
\end{assum}

Under Assumption \ref{assum:volume}, 
the diameter of $M$ is infinity.
We see from \cite[Proposition 5.2]{GH14} that 
if $M$ is non-compact and connected, 
and if $B(x,r)$ is relatively compact for any $x\in M$ and $r>0$, 
then Assumption \ref{assum:volume} is fulfilled 
under the condition that for some $c>0$ and $d>0$, 
$$\frac{V(x,R)}{V(x,r)}\leq c\left(\frac{R}{r}\right)^d \quad  \text{for all $x\in M$ and  $0<r<R$}.$$

\begin{assum}\label{assum:heat-kernel}
There exist a properly exceptional Borel set $N$ and a nonnegative symmetric kernel $p(t,x,y)$ 
on $(0,\infty)\times (M\setminus N) \times (M\setminus N)$ such that 
$$P_x(X_t\in A)=\int_A p(t,x,y)\,m({\rm d}y), \quad t\geq 0, \ A\in {\cal B}(M)$$
for any $x\in M\setminus N$ and 
$$p(t+s,x,y)=\int_M p(t,x,z)p(s,z,y)\,m({\rm d}z)$$
for any $x,y\in M\setminus N$ and $t,s>0$.
\end{assum}

The nonnegative symmetric kernel $p(t,x,y)$ in Assumption \ref{assum:heat-kernel} 
is called the {\it heat kernel} of ${\bf M}$. 

\begin{assum}\label{assum:estimate}
There exist positive constants $L_1$ and $L_2$ such that 
for any $(t,x,y)\in (0,\infty)\times (M\setminus N)\times (M\setminus N)$,
\begin{equation}\label{estimate}
\begin{split}
&L_1\min\left\{\frac{1}{V(x,\phi^{-1}(t))}, \frac{t}{V(x,d(x,y))\phi(d(x,y))}\right\}\\
&\leq 
p(t,x,y) 
\leq L_2\min\left\{\frac{1}{V(x,\phi^{-1}(t))}, \frac{t}{V(x,d(x,y))\phi(d(x,y))}\right\}.
\end{split}
\end{equation}
Here $\phi(r)$ is some positive increasing function on $[0,\infty)$ such that 
\begin{equation}\label{eq:scale-growth}
\phi(0)=0, \quad c_3\left(\frac{R}{r}\right)^{d_3}\leq \frac{\phi(R)}{\phi(r)}\leq c_4\left(\frac{R}{r}\right)^{d_4} \quad (0<r<R)
\end{equation}
for some constants  $c_3\in (0,1]$, $c_4\in [1,\infty)$, $d_3>0$ and $d_4\geq d_3$. 
\end{assum}

It follows by (\ref{eq:scale-growth}) that 
\begin{equation}\label{eq:scale-growth-inverse}
\frac{1}{c_4^{1/d_4}}\left(\frac{R}{r}\right)^{1/d_4}
\leq \frac{\phi^{-1}(R)}{\phi^{-1}(r)}
\leq\frac{1}{c_3^{1/d_3}}\left(\frac{R}{r}\right)^{1/d_3} \ (0<r<R).
\end{equation}

A positive function $r(t)$ on $(0,\infty)$ is called a \textit{lower rate function} of ${\bf M}$ if 
$$P_x\left(\text{there exists $T>0$ such that $d(X_t,x)>r(t)$ for all $t>T$}\right)=1$$
for q.e.\  $x\in M$. 
For a lower rate  function $r(t)$, the probability 
$$q_r(t,x):=P_x\left(\text{$d(X_u,x)\leq r(u)$ for some $u>t$}\right) \quad (x\in M, \ t>0),$$ 
tends to $0$ as $t\rightarrow\infty$ for q.e.\ $x\in M$.
In what follows, we find the decay rate of $q_r(t,x)$ as $t\rightarrow \infty$.

\subsection{Transient case}
We first assume that $d_1>d_4$. 
Then ${\bf M}$ is transient in the sense that 
$$P_x\left(\lim_{t\rightarrow\infty}X_t=\Delta\right)=1\quad \text{for q.e.\ $x\in M$}.$$
If $B(x,r)$ is relatively compact for any $x\in M$ and $r>0$, 
then the transience is equivalent to 
$$P_x\left(\lim_{t\rightarrow\infty}d(x,X_t)=\infty\right)=1\quad \text{for q.e.\ $x\in M$}.$$
If we assume in addition that $M$ is connected and  
$$V(x,r)\asymp V(r)$$ 
for some positive function $V(r)$ on $(0,\infty)$,
then the following integral test on the lower rate functions 
is derived in \cite{SW15}:
let $g(t)$ be a positive decreasing function on $(0,\infty)$ 
such that $g(t)\rightarrow 0$ as $t\rightarrow\infty$ and $\varphi(t)=\phi^{-1}(t)g(t)$. 
Then 
\begin{equation}\label{0-1}
P_x\left(\text{there exists $T>0$ such that $d(X_u,x)>\varphi(u)$ for all $u\geq T$}\right)=\text{$1$ or $0$}
\end{equation}
for q.e.\ $x\in M$ according as 
\begin{equation}\label{eq:1-prob}
\int_1^{\infty}\frac{V(\varphi(t))}{\phi(\varphi(t))}\frac{{\rm d}t}{V(\phi^{-1}(t))}<\infty \ \text{or $=\infty$}.
\end{equation}

Fix a positive decreasing function $g(t)$ on $(0,\infty)$ such that 
$g(t)\rightarrow 0$ as $t\rightarrow\infty$.  
Define for $t>0$ and $c>1$,
$$
R_{c,t}=\inf\left\{\frac{g(u)}{g(v)} \mid 1\leq \frac{u}{v}\leq c, \ v\geq t\right\}.
$$
We then have

\begin{thm}\label{asymp-thm}
Suppose that Assumptions {\rm \ref{assum:volume}-- \ref{assum:estimate}} are fulfilled and $d_1>d_4$.  
Let $g(t)$ be a positive decreasing function on $(0,\infty)$ 
such that $g(t)\rightarrow 0$ as $t\rightarrow\infty$ and $\varphi(t)=\phi^{-1}(t)g(t)$. 
For q.e.\ $x\in M$, if 
$$\int_1^{\infty}\frac{V(x,\varphi(t))}{\phi(\varphi(t))}\frac{{\rm d}t}{V(x,\phi^{-1}(t))}<\infty$$
 and $\lim_{c\rightarrow 1+0}\left(\lim_{t\rightarrow\infty}R_{c,t}\right)=1$,  then 
$$q_{\varphi}(t,x)\asymp \int_t^{\infty}\frac{V(x,\varphi(s))}{\phi(\varphi(s))}\frac{{\rm d}s}{V(x,\phi^{-1}(s))} 
\quad (t\rightarrow\infty).$$
\end{thm}

\begin{rem}\rm \ Under Assumptions \ref{assum:volume}--\ref{assum:estimate},  
$({\cal E}, {\cal F})$ is conservative and 
the heat kernel $p(t,x,y)$ satisfies (\ref{estimate}) for all $(t,x,y)\in (0,\infty)\times M\times M$ 
by \cite[Theorem 1.13 and Proposition 3.1]{CKW16}. 
We can thus take a version of the process ${\bf M}$ 
such that Theorem \ref{asymp-thm} is valid for any $x\in M$.
Moreover, there exists a nonnegative measurable function $J(x,y)$ on $M\times M\setminus {\rm diag}$ such that  
$J({\rm d}x{\rm d}y)=J(x,y)\,m({\rm d}x)m({\rm d}y)$ and
$$J(x,y)\asymp \frac{1}{V(x, d(x,y))\phi(d(x,y))}.$$
\end{rem}

\subsection{Recurrent case}
We next assume that for some positive constants $c_{v,1}$ and $c_{v,2}$, 
\begin{equation}\label{eq:n-p-r}
c_{v,1}\phi(r)\leq V(x,r)\leq c_{v,2}\phi(r) \quad \text{for all $x\in M$ and $r>0$}.
\end{equation}
Then ${\bf M}$ is irreducible recurrent by \cite[Remark 2.2]{SW15}. Moreover, we can show that ${\bf M}$ can not hit any point by following the proof of \cite[Proposition 4.24]{SW15} and using Lemma \ref{prop:hit-prob} below. These two properties imply that 
$$P_x\left(\text{$\liminf_{t\rightarrow\infty}d(x,X_t)=0$ and  $d(x,X_t)>0$ for all $t>0$}\right)=1 
\quad \text{for q.e.\ $x\in M$}.$$ 
We proved in \cite{SW15} that if $M$ is, in addition, connected, then (\ref{0-1}) is valid according as 
\begin{equation}\label{eq:1-prob-critical}
\int_1^{\infty}\frac{{\rm d}t}{t|\log g(t)|}<\infty \ \text{or $=\infty$}.
\end{equation}

For $t>0$ and $c>1$, we define 
$$R_{c,t}=
\sup\left\{\frac{|\log g(u)|}{|\log g(v)|} \mid 1\leq \frac{u}{v}\leq c, \ v\geq t\right\}.
$$

\begin{thm}\label{asymp-thm-critical}
Suppose that Assumptions {\rm \ref{assum:volume}-- \ref{assum:estimate}} and {\rm (\ref{eq:n-p-r})} 
are fulfilled. 
Let $g(t)$ be a positive decreasing function on $(0,\infty)$ 
such that $g(t)\rightarrow 0$ as $t\rightarrow\infty$ and $\varphi(t)=\phi^{-1}(t)g(t)$. 
If the integral in {\rm (\ref{eq:1-prob-critical})} is convergent  and 
$\lim_{c\rightarrow 1+0}(\lim_{t\rightarrow\infty}R_{c,t})=1$,  
then 
$$q_{\varphi}(t,x)\asymp  \int_t^{\infty}\frac{{\rm d}s}{s|\log g(s)|} \quad (t\rightarrow\infty)$$
for q.e.\  $x\in M$. 
\end{thm}

Remarks similar to those just after Theorem \ref{asymp-thm} 
are also valid in  Theorem \ref{asymp-thm-critical}.

\section{Examples}\label{sect:exam}

We first apply Theorems \ref{asymp-thm} and \ref{asymp-thm-critical} 
to some lower rate functions. 

\begin{exam}\label{exam:order}\rm 
Suppose  that for some positive constants $\alpha_1$ and $\alpha_2$, 
\begin{equation}\label{ex-volume}
V(x,r)\asymp r^{\alpha_1}{\bf 1}_{\{r<1\}}+r^{\alpha_2}{\bf 1}_{\{r\geq 1\}}.
\end{equation}
Then Assumption \ref{assum:volume} is fulfilled 
by $d_1=\alpha_1\wedge \alpha_2$ and $d_2=\alpha_1\vee \alpha_2$. 
We also suppose that 
Assumptions \ref{assum:heat-kernel} and \ref{assum:estimate} are 
satisfied by the functions 
\begin{equation}\label{ex-scale}
\phi(r)=r^{\beta_1}{\bf 1}_{\{r<1\}}+r^{\beta_2}{\bf 1}_{\{r\geq 1\}}
\end{equation}
for some positive constants $\beta_1$ and $\beta_2$.  
Then $d_3=\beta_1\wedge \beta_2$ and $d_4=\beta_1\vee \beta_2$. 
Let  $g(t)$ be a positive decreasing function on $(0,\infty)$ 
such that $g(t)\rightarrow 0$ as $t\rightarrow\infty$ and $\varphi(t)=\phi^{-1}(t)g(t)$.

We first assume that  
$\alpha_1\wedge \alpha_2>\beta_1\vee \beta_2$.  
If  $g(t)$ satisfies the full assumptions in Theorem \ref{asymp-thm},  then as $t\rightarrow\infty$,
$$
q_{\varphi}(t,x)
\asymp 
\begin{cases}
\displaystyle \int_t^{\infty}\varphi(s)^{\alpha_2-\beta_2}\frac{{\rm d}s}{s^{\frac{\alpha_2}{\beta_2}}} 
& \text{if $\varphi(t)\rightarrow\infty$ \ ($t\rightarrow\infty$),}\bigskip \\
\displaystyle \frac{1}{t^{\frac{\alpha_2-\beta_2}{\beta_2}}} 
&  \text{if $\varphi(t)\asymp 1$ \ ($t\rightarrow\infty$),} \bigskip \\
\displaystyle \int_t^{\infty}\varphi(s)^{\alpha_1-\beta_1}\frac{{\rm d}s}{s^{\frac{\alpha_2}{\beta_2}}} 
& \text{if $\varphi(t)\rightarrow 0$ \ ($t\rightarrow\infty$).}
\end{cases}$$
This implies that if $\varphi(t)=t^{1/\beta_2}/(\log t)^{\frac{1+\varepsilon}{\alpha_2-\beta_2}}$ 
for some $\varepsilon>0$, then 
\begin{equation}\label{eq:asymp-ex}
q_{\varphi}(t,x)
\asymp \frac{1}{\varepsilon(\log t)^{\varepsilon}} \quad (t\rightarrow\infty).
\end{equation}
On the other hand, if  $\varphi(t)=t^p$ for some $p<1/\beta_2$, then 
$$
q_{\varphi}(t,x)
\asymp 
\begin{cases}
\displaystyle \frac{1}{t^{\left(\frac{1}{\beta_2}-p\right)(\alpha_2-\beta_2)}}
& \text{if $\displaystyle 0\leq p< \frac{1}{\beta_2}$,}\bigskip \\
\displaystyle \frac{1}{t^{\frac{1}{\beta_2}(\alpha_2-\beta_2)-p(\alpha_1-\beta_1)}}
& \text{if $p<0$.}
\end{cases}$$

We next assume that $\alpha_1=\beta_1$ and $\alpha_2=\beta_2$. 
If $g(t)$ satisfies the full assumptions 
in Theorem \ref{asymp-thm-critical}, then 
$$
q_{\varphi}(t,x)\asymp\int_t^{\infty}\frac{{\rm d}s}{s|\log g(s)|} \quad (t\rightarrow\infty).$$  
Hence if $\varphi(t)=e^{-(\log t)^{1+\varepsilon}}$ for some $\varepsilon>0$, then 
we get (\ref{eq:asymp-ex}).
On the other hand, if $\varphi(t)=e^{-t^p}$ for some $p>0$, then 
$$
q_{\varphi}(t,x)
\asymp \frac{1}{t^{p}} \quad (t\rightarrow\infty).
$$

We now provide examples satisfying  
Assumptions \ref{assum:heat-kernel} and \ref{assum:estimate}. 
Suppose that 
\begin{itemize}
\item $B(x,r)$ is relatively compact for any $x\in M$ and $r>0$;
\item the distance $d$ on $M$ is geodesic: for any $x,y\in M$, 
there exists a continuous map $\gamma: [0,1]\rightarrow M$ such that 
$\gamma(0)=x$, $\gamma(1)=y$ and 
$$d(\gamma(t),\gamma(s))=|t-s|d(x,y) \quad \text{for any $t,s\in [0,1]$}.$$
\end{itemize}
Let $V(r)$ and $\phi(r)$ satisfy (\ref{ex-volume}) and (\ref{ex-scale}) with $\beta_1, \beta_2\in (0,2)$, 
respectively.  
Let $c(x,y)$ be a  uniformly positive and bounded function on $M\times M$
and 
$$
{\cal E}(u,v)=\iint_{M\times M\setminus {\rm diag}}\frac{(u(x)-u(y))(v(x)-v(y))}{V(d(x,y))\phi(d(x,y))}c(x,y)\,m({\rm d}x)m({\rm d}y)
$$
for $u,v\in L^2(M;m)$ if the right hand side above makes sense. 
Denote by ${\cal F}$ the ${\cal E}_1$-closure of 
the totality of Lipschitz continuous functions on $M$ with compact support. 
Then $({\cal E}, {\cal F})$ is a regular Dirichlet form on $L^2(M;m)$ 
so that there exists an associated  $m$-symmetric Hunt process ${\bf M}$ on $M$ of pure jump type. 
Furthermore, Assumptions \ref{assum:heat-kernel} and \ref{assum:estimate} are fulfilled 
according to \cite{CK08} (see also \cite{C09, CK03, KKW15}). 
\end{exam}

\begin{exam}\label{ex:sub-diff} \rm 
Suppose that Assumption \ref{assum:volume} is fulfilled. 
Let ${\bf M}$ be a conservative $m$-symmetric diffusion process on $M$ such that 
the associated Dirichlet form $({\cal E}, {\cal F})$ is regular on $L^2(M;m)$. 
Suppose further that ${\bf M}$ admits the heat kernel $p(t,x,y)$ 
such that for some positive constants $c_i$, $C_i$  $(1\leq i\leq 4)$ 
and $\beta_1, \beta_2\in [2,\infty)$, the following hold:
\begin{itemize}
\item if $0<t\leq 1\vee d(x,y)$, then  
\begin{equation*}
\begin{split}
\frac{c_1}{V(x,t^{1/\beta_1})}
\exp\left\{-C_1\left(\frac{d(x,y)^{\beta_1}}{t}\right)^{\frac{1}{\beta_1-1}}\right\}
&\leq p(t,x,y)\\
&\leq \frac{c_2}{V(x,t^{1/\beta_1})}
\exp\left\{-C_2\left(\frac{d(x,y)^{\beta_1}}{t}\right)^{\frac{1}{\beta_1-1}}\right\};
\end{split}
\end{equation*}
\item if $t\geq 1\vee d(x,y)$, then  
\begin{equation*}
\begin{split}
\frac{c_3}{V(x,t^{1/\beta_2})}
\exp\left\{-C_3\left(\frac{d(x,y)^{\beta_2}}{t}\right)^{\frac{1}{\beta_2-1}}\right\}
&\leq p(t,x,y)\\
&\leq \frac{c_4}{V(x,t^{1/\beta_2})}
\exp\left\{-C_4\left(\frac{d(x,y)^{\beta_2}}{t}\right)^{\frac{1}{\beta_2-1}}\right\}.
\end{split}
\end{equation*}
\end{itemize}
If $\beta_1=\beta_2=2$, then the heat kernel $p(t,x,y)$ admits the so-called Gaussian estimates. 
Here we allow $\beta_1$ and  $\beta_2$ to be different and greater than 2. 
Such situation occurs for a class of symmetric diffusion processes on (unbounded) fractals and fractal-like spaces such as Sierpi\'nski carpets 
and pre-carpets
(see, e.g., \cite{AB15, BB99, BB00, BBK06, O96}. See also the summary  just after  \cite[Proposition 5.5]{K13} 
on the history of the analysis on Sierpi\'nski carpets). 

For $\gamma\in (0,1)$, we let ${\bf M}^{(1)}$ be 
a $\gamma$-stable subordinated diffusion process of ${\bf M}$ 
(see Appendix \ref{sect:sub} below for definition). 
According to \cite[Theorem 2.1]{O02}, 
the associated Dirichlet form $({\cal E}^{(1)}, {\cal F}^{(1)})$ is regular on $L^2(M;m)$ and 
\begin{equation}\label{eq:dirichlet-sub}
{\cal E}^{(1)}(u,u)
\asymp \iint_{M\times M\setminus{\rm diag}}\frac{(u(x)-u(y))^2}{V(x,d(x,y))\phi(d(x,y))}\,m({\rm d}x)m({\rm d}y)
\end{equation}
for any $u\in {\cal F}^{(1)}\cap C_0(M)$, where 
$$\phi(r)=r^{\gamma \beta_1}{\bf 1}_{\{r<1\}}+r^{\gamma\beta_2}{\bf 1}_{\{r\geq 1\}}.$$
Furthermore, ${\bf M}^{(1)}$ admits the heat kernel $q(t,x,y)$ such that 
\begin{equation}\label{eq:heat-sub}
q(t,x,y)\asymp \min\left\{\frac{1}{V(x,\phi^{-1}(t))},\frac{t}{V(x,d(x,y))\phi(d(x,y))}\right\}.
\end{equation}
Therefore, Assumptions \ref{assum:heat-kernel} and \ref{assum:estimate} are valid 
for ${\bf M}^{(1)}$.  
We note that if $\alpha_1=\alpha_2$ and $\beta_1=\beta_2$, 
then $q(t,x,y)$ is already computed in \cite{BSS03} and \cite{K03}. 
We show (\ref{eq:dirichlet-sub}) and  (\ref{eq:heat-sub}) in Appendix \ref{sect:sub} below. 
\end{exam}

We finally apply Theorem \ref{asymp-thm} 
to subordinated diffusion processes under the non-uniform volume growth condition.
\begin{exam} \rm 
Suppose  that $M={\mathbb R}^d$ and $m$ is the Lebesgue measure on ${\mathbb R}^d$ 
(${\rm d}x$ in notation). 
Let $h(x)$ be a positive function on ${\mathbb R}^d$ 
such that $h(x)\asymp (1+|x|^2)^{\alpha/2}$ for some $\alpha>-d/2$ 
and $\mu({\rm d}x)=h(x)^2\,{\rm d}x$. 
We denote by $C_0^{\infty}({\mathbb R}^d)$ the totality of smooth functions on ${\mathbb R}^d$ with compact support. 
Let ${\cal E}$ be a bilinear form on 
$C_0^{\infty}({\mathbb R}^d)\times C_0^{\infty}({\mathbb R}^d)$ 
defined by 
$${\cal E}(u,v)=\frac{1}{2}\int_{{\mathbb R}^d} 
\nabla u(x)\cdot \nabla v(x)\, \mu({\rm d}x), \quad  u,v\in C_0^{\infty}({\mathbb R}^d)$$ 
and let ${\cal F}$ be the closure of $C_0^{\infty}({\mathbb R}^d)$ with respect to the norm 
$$\|u\|_{{\cal E}_1}=\sqrt{{\cal E}(u,u)+\|u\|_{L^2({\mathbb R}^d;\mu)}^2}.$$
Then $({\cal E}, {\cal F})$ is a regular Dirichlet form on $L^2({\mathbb R}^d;\mu)$ 
such that there exists a $\mu$-symmetric diffusion process 
${\bf M}$ on ${\mathbb R}^d$. 
According to \cite[Subsection 4.3]{GS05} and \cite[Section 4]{St96} 
(see also \cite[Corollary 6.11]{G06}), 
we have 
$$V(x,r)=\mu(B(x,r))\asymp r^d(1+r+|x|)^{2\alpha}$$ 
so that Assumption \ref{assum:volume} is valid for $d_1=d$ and $d_2=d+2\alpha$. 
Furthermore,  the associated heat kernel $p(t,x,y)$ satisfies for any $x,y\in {\mathbb R}^d$ and $t>0$,
\begin{equation*}
\begin{split}
\frac{c_1}{t^{d/2}(1+\sqrt{t}+|x|)^{2\alpha}}\exp\left(-C_1\frac{d(x,y)^2}{t}\right)
&\leq p(t,x,y)\\
&\leq \frac{c_2}{t^{d/2}(1+\sqrt{t}+|x|)^{2\alpha}}\exp\left(-C_2\frac{d(x,y)^2}{t}\right).
\end{split}
\end{equation*}

For $\gamma\in (0,1)$, we let ${\bf M}^{(1)}$ be 
a $\gamma$-stable subordinated diffusion process of ${\bf M}$. 
Then by the same argument as in Example \ref{ex:sub-diff}, 
the associated Dirichlet form $({\cal E}^{(1)}, {\cal F}^{(1)})$ is regular on $L^2({\mathbb R}^d;\mu)$ and 
$$
{\cal E}^{(1)}(u,u)
\asymp \iint_{{\mathbb R}^d\times {\mathbb R}^d\setminus{\rm diag}}
\frac{(u(x)-u(y))^2}{|x-y|^{d+2\gamma}(1+|x-y|+|x|)^{2\alpha}}\,\mu({\rm d}x)\mu({\rm d}y)
$$
for any $u\in C_0^{\infty}({\mathbb R}^d)$.
Furthermore, ${\bf M}^{(1)}$ admits the heat kernel $q(t,x,y)$ such that 
$$
q(t,x,y)\asymp 
\min\left\{\frac{1}{t^{d/(2\gamma)}(1+t^{1/(2\gamma)}+|x|)^{2\alpha}},
\frac{t}{(1+|x-y|+|x|)^{2\alpha}|x-y|^{d+2\gamma}}\right\}.
$$
Therefore, ${\bf M}^{(1)}$ is transient for $d+2\alpha>2\gamma$. 

Let  $g(t)$ be a positive decreasing function on $(0,\infty)$ 
such that $g(t)\rightarrow 0$ as $t\rightarrow\infty$ and $\varphi(t)=t^{1/(2\gamma)} g(t)$. 
If $d>2\gamma$ and $d+2\alpha>2\gamma$, then under the full conditions  in Theorem \ref{asymp-thm},  
$$q_{\varphi}(t,x)\asymp \int_t^{\infty}\frac{\varphi(s)^{d-2\gamma}}{s^{d/(2\gamma)}}
\left(\frac{1+\varphi(s)+|x|}{1+s^{1/(2\gamma)}+|x|}\right)^{2\alpha}\,{\rm d}s \quad (t\rightarrow\infty).$$
This implies that if $\varphi(t)=t^{1/(2\gamma)}/(\log t)^{\frac{1+\varepsilon}{d+2\alpha-2\gamma}}$ 
for some $\varepsilon>0$, then 
\begin{equation}\label{eq:asymp-ex}
q_{\varphi}(t,x)
\asymp \frac{1}{\varepsilon(\log t)^{\varepsilon}} \quad (t\rightarrow\infty).
\end{equation}
On the other hand, if  $\varphi(t)=t^p$ for some $p<1/(2\gamma)$, then 
$$
q_{\varphi}(t,x)
\asymp 
\begin{cases}
\displaystyle \frac{1}{t^{\left(\frac{1}{2\gamma}-p\right)(d+2\alpha-2\gamma)}},
& \text{if $\displaystyle 0\leq p< \frac{1}{2\gamma}$,}\bigskip \\
\displaystyle \frac{1}{t^{\frac{1}{2\gamma}(d+2\alpha-2\gamma)-p(d-2\gamma)}},
& \text{if $p<0$.}
\end{cases}$$
\end{exam}

\section{Hitting time distributions}\label{sect:hit}
Throughout this section, we assume 
the full conditions in Theorem \ref{asymp-thm}. 
For simplicity, we also assume that $N=\emptyset$ in Assumptions \ref{assum:heat-kernel} and \ref{assum:estimate}. 
For the proof of  Theorem \ref{asymp-thm}, 
we first give estimates of the hitting time distributions 
to closed balls during finite time interval. 
To do so, we use the next lemma which goes back to \cite{K97}.

\begin{lem}\label{lem:hit} {\rm (\cite[Lemma 4.19]{SW15})}.
Let $a$, $b$, $c$ and $r$ be positive constants. 
Then for any $x\in M$,
\begin{equation*}
\begin{split}
\frac{\displaystyle \int_a^b P_x(d(X_u,x)\leq r)\,{\rm d}u}
{\displaystyle 2\int_{0}^{b-a}\sup_{d(y,x)\leq r}P_y(d(X_u,y)\leq 2r)\,{\rm d}u}
&\leq 
P_x\left(\text{$d(X_s,x)\leq r$ for some $s\in (a,b]$}\right)\\
&\leq \frac{\displaystyle \int_a^{b+c}P_x(d(X_u,x)\leq 2r)\,{\rm d}u}
{\displaystyle \int_0^c \inf_{d(y,x)\leq r}P_y(d(X_u,y)\leq r)\,{\rm d}u}.
\end{split}
\end{equation*}
\end{lem}

By the same way as in \cite[Lemma 4.20]{SW15}, we have 
\begin{equation}\label{eq:ball-inside}
L_1\min\left\{1,\frac{V(x,r)}{V(x,\phi^{-1}(t))}\right\}\leq P_x(d(X_t,x)\leq r)
\leq L_2\min\left\{1,\frac{V(x,r)}{V(x,\phi^{-1}(t))}\right\}
\end{equation}
for any $x\in M$, $r>0$ and $t>0$. 
Using this inequality, we have the following two lemmas.

\begin{lem}\label{prop:hit-prob}
Let $a$, $b$, $c$ and $r$ be positive constants. 
If $\phi(r)\leq a\wedge c$, then for any $x\in M$,
$$P_x\left(\text{$d(X_s,x)\leq r$ for some $s\in (a,b]$}\right)
\leq K_1 \frac{V(x,r)}{\phi(r)}
\int_a^{b+c}\frac{{\rm d}u}{V(x,\phi^{-1}(u))},$$
where $K_1=2^{d_2}c_2L_2/L_1$.
\end{lem}

\pf \ Suppose that $\phi(r)\leq a\wedge c$. 
Then by (\ref{eq:ball-inside}) and (\ref{eq:volume-growth}), we have 
\begin{equation}\label{eq:hit-prob-est-1}
\begin{split}
&\int_a^{b+c}P_x(d(X_u,x)\leq 2r)\,{\rm d}u
\leq L_2\int_a^{b+c}\min\left\{1,\frac{V(x,2r)}{V(x,\phi^{-1}(u))}\right\}\,{\rm d}u\\
&\leq 2^{d_2}c_2L_2\int_a^{b+c}\min\left\{1,\frac{V(x,r)}{V(x,\phi^{-1}(u))}\right\}\,{\rm d}u
=2^{d_2}c_2L_2V(x,r)\int_a^{b+c}\frac{{\rm d}u}{V(x,\phi^{-1}(u))}.
\end{split}
\end{equation}
We also see by (\ref{eq:ball-inside}) that if $0<u\leq \phi(r)$, then for any $y\in M$, 
$$P_y(d(X_u,y)\leq r)\geq L_1\min\left\{1,\frac{V(y,r)}{V(y,\phi^{-1}(u))}\right\}=L_1,$$
which implies that 
$$\int_0^c\inf_{d(y,x)\leq r}P_y(d(X_u,y)\leq r)\,{\rm d}u
\geq \int_0^{\phi(r)}\inf_{d(y,x)\leq r}P_y(d(X_u,y)\leq r)\,{\rm d}u\geq L_1\phi(r).$$
Therefore, the proof is completed by Lemma \ref{lem:hit}. 
\qed

\begin{lem}\label{lem:hit-lower}
Let $a$, $b$ and $r$ be positive constants. If $\phi(r)\leq a$  and $\phi(2r)\leq b-a$, then for any $x\in M$,
\begin{equation}\label{eq:est-hit-lower-2}
P_x\left(\text{$d(X_s,x)\leq r$ for some $s\in (a,b]$}\right)
\geq K_2 \frac{V(x,r)}{\phi(r)}\int_a^b\frac{{\rm d}u}{V(x,\phi^{-1}(u))},
\end{equation}
where 
$$K_*=\frac{c_4^{d_1/d_4}}{c_1}\frac{d_4}{d_1-d_4}, \quad 
K_2=\frac{L_1}{L_2}\frac{1}{c_42^{d_4+1}(1+(c_2)^23^{d_2}K_*)}.$$
\end{lem}

\pf \ 
Suppose that $\phi(r)\leq a$  and $\phi(2r)\leq b-a$. 
Then by (\ref{eq:ball-inside}), 
\begin{equation*}
\begin{split}
\int_a^b P_x(d(X_u,x)\leq r)\,{\rm d}u
&\geq L_1\int_a^b \min\left\{1,\frac{V(x,r)}{V(x,\phi^{-1}(u))}\right\}\,{\rm d}u\\
&=L_1V(x,r)\int_a^b\frac{{\rm d}u}{V(x,\phi^{-1}(u))}.
\end{split}
\end{equation*}
If $d(y,x)\leq r$, then we have $B(y,2r)\subset B(x,3r)$ by the triangle inequality so that $V(y,2r)\leq  V(x,3r)$.  
If we assume in addition that $u\geq \phi(2r)$, then 
$$\phi^{-1}(u)-d(y,x)\geq \phi^{-1}(u)-r\geq\frac{1}{2}\phi^{-1}(u)$$
so that 
$$V(y,\phi^{-1}(u))\geq V(x,\phi^{-1}(u)-d(y,x))
\geq V(x,\phi^{-1}(u)/2)\geq \frac{1}{c_22^{d_2}}V(x,\phi^{-1}(u))$$
by 
the triangle inequality and (\ref{eq:volume-growth}). 
Hence (\ref{eq:ball-inside}) implies that 
$$P_y(d(X_u,y)\leq 2r)\leq L_2\frac{V(y,2r)}{V(y,\phi^{-1}(u))}\leq L_2c_22^{d_2}\frac{V(x,3r)}{V(x,\phi^{-1}(u))}.$$
Since this inequality yields that 
$$
\int_{0}^{b-a}\sup_{d(y,x)\leq r}P_y(d(X_u,y)\leq 2r)\,{\rm d}u
\leq L_2\left(\phi(2r)+c_22^{d_2}V(x,3r)\int_{\phi(2r)}^{b-a} \frac{{\rm d}u}{V(x,\phi^{-1}(u))}\right),
$$
we have by Lemma \ref{lem:hit}, 
\begin{equation}\label{eq:est-hit-lower-1}
\begin{split}
&P_x\left(\text{$d(X_s,x)\leq r$ for some $s\in (a,b]$}\right)\\
&\geq \frac{L_1}{2L_2}V(x,r)\int_a^b\frac{{\rm d}u}{V(x,\phi^{-1}(u))}
\frac{1}{\displaystyle \phi(2r)+c_22^{d_2}V(x,3r)\int_{\phi(2r)}^{b-a} \frac{{\rm d}u}{V(x,\phi^{-1}(u))}}.
\end{split}
\end{equation}

Fix $\theta>1$ and take $t_m=t\theta^m \ (m\geq 0)$. Then 
\begin{equation}\label{eq:int-sum}
\begin{split}
\int_t^{\infty}\frac{{\rm d}u}{V(x,\phi^{-1}(u))}
=\sum_{m=0}^{\infty}\int_{t_m}^{t_{m+1}}\frac{{\rm d}u}{V(x,\phi^{-1}(u))}
&\leq \sum_{m=0}^{\infty}\frac{t_{m+1}-t_m}{V(x,\phi^{-1}(t_m))}\\
&=t\sum_{m=0}^{\infty}\frac{\theta^{m+1}-\theta^{m}}{V(x,\phi^{-1}(t_m))}.
\end{split}\end{equation}
Since $d_1>d_4$ by assumption and  
$$\frac{V(x,\phi^{-1}(t_m))}{V(x,\phi^{-1}(t))}\geq \frac{c_1}{c_4^{d_1/d_4}}(\theta^m)^{d_1/d_4}$$
by (\ref{eq:volume-growth}) and (\ref{eq:scale-growth-inverse}),
the last expression of (\ref{eq:int-sum}) is less than 
$$
\frac{c_4^{d_1/d_4}}{c_1}\frac{t}{V(x,\phi^{-1}(t))}(\theta-1)\sum_{m=0}^{\infty}(\theta^m)^{1-d_1/d_4}
=\frac{c_4^{d_1/d_4}}{c_1}\frac{\theta-1}{1-\theta^{1-d_1/d_4}}\frac{t}{V(x,\phi^{-1}(t))}.
$$
Then by letting $\theta\rightarrow1+0$, we get
$$\int_t^{\infty}\frac{{\rm d}u}{V(x,\phi^{-1}(u))}\leq K_*\frac{t}{V(x,\phi^{-1}(t))}$$
so that 
$$\int_{\phi(2r)}^{b-a}\frac{{\rm d}u}{V(x,\phi^{-1}(u))}
\leq \int_{\phi(2r)}^{\infty}\frac{{\rm d}u}{V(x,\phi^{-1}(u))}
\leq K_* \frac{\phi(2r)}{V(x,2r)}.$$
Since  (\ref{eq:volume-growth}) and (\ref{eq:scale-growth}) imply that 
$$\frac{V(x,3r)}{V(x,2r)}\leq c_2\left(\frac{3}{2}\right)^{d_2}$$
and $\phi(2r)\leq c_42^{d_4}\phi(r)$, respectively, 
we obtain (\ref{eq:est-hit-lower-2}) by  (\ref{eq:est-hit-lower-1}).
\qed
\medskip

\section{Proof of Theorem \ref{asymp-thm}}\label{sect:proof}
In this section, we prove Theorem \ref{asymp-thm} 
by using the results in Section \ref{sect:hit}. 
More precisely, we show that 
\begin{equation}\label{eq:sup}
\limsup_{t\rightarrow\infty}
\frac{q_{\varphi}(t,x)}
{\displaystyle \int_t^{\infty}\frac{V(x,\varphi(s))}{\phi(\varphi(s))}\frac{{\rm d}s}{V(x,\phi^{-1}(s))}}
\leq \frac{L_2}{L_1}\frac{2^{d_2}(c_2)^2}{c_3^{d_2/d_3}}
\end{equation}
and 
\begin{equation}\label{eq:inf}
\liminf_{t\rightarrow\infty}
\frac{q_{\varphi}(t,x)}
{\displaystyle \int_t^{\infty}\frac{V(x,\varphi(s))}{\phi(\varphi(s))}\frac{{\rm d}s}{V(x,\phi^{-1}(s))}}
\geq \frac{L_1}{L_2}\frac{(c_1)^2(c_3)^{3d_2/d_3-1}}
{2^{d_4+1}(c_2c_4)^2}\frac{d_1-d_4}{(d_1-d_4)c_1+3^{d_2}d_4(c_2)^2c_4^{d_1/d_4}}.
\end{equation}
Throughout this section, we keep the same setting as in Section \ref{sect:hit}. 

\subsection{Proof of (\ref{eq:sup})}
For fixed constants $t>0$ and $c\in (1,2)$, 
we define a sequence $\{n_k\}_{k=0}^{\infty}$ by $n_k=tc^k \ (k\geq 0)$.
Then 
\begin{equation}\label{eq:sup-est-1}
\begin{split}
q_{\varphi}(t,x)
&=P_x\left(\text{$d(X_u,x)\leq \varphi(u)$ for some $u>t$}\right)\\
&=P_x\left(\bigcup_{k=0}^{\infty}\left\{\text{$d(X_u,x)\leq \varphi(u)$ for some $u\in (n_k,n_{k+1}]$}\right\}\right)\\
&\leq\sum_{k=0}^{\infty}P_x\left(\text{$d(X_u,x)\leq \varphi(u)$ for some $u\in (n_k,n_{k+1}]$}\right).
\end{split}
\end{equation}
To obtain an upper bound of the last term of (\ref{eq:sup-est-1}), we show

\begin{lem}\label{lem:hit-prob-1}
For each $c\in (1,2)$, there exists $T_c>0$ such that for all $t\geq T_c$, 
\begin{equation}\label{eq:hit-prob-1}
\begin{split}
&P_x\left(\text{$d(X_u,x)\leq \varphi(u)$ for some $u\in (n_k,n_{k+1}]$}\right)\\
&\leq  \frac{K_1c_2}{c_3^{d_2/d_3}} \frac{c^{1+d_2/d_3}}{(R_{c,t})^{d_2}}
\int_{n_k}^{n_{k+1}}\frac{V(x,\varphi(u))}{\phi(\varphi(u))}\frac{{\rm d}u}{V(x,\phi^{-1}(u))} 
\end{split}
\end{equation}
for any $x\in M$ and $k\geq 0$.
\end{lem}

\pf \ For any $u\in (n_k, n_{k+1}]$, we have by (\ref{eq:scale-growth-inverse}), 
\begin{equation}\label{eq:comp-radius}
\begin{split}
\varphi(u)=\phi^{-1}(u)g(u)
\leq \left(\frac{c}{c_3}\right)^{1/d_3}\phi^{-1}(n_k)g(n_k)
=\left(\frac{c}{c_3}\right)^{1/d_3}\varphi(n_k),
\end{split}
\end{equation}
which implies that 
\begin{equation}\label{eq:hit-prob-0}
\begin{split}
&P_x\left(\text{$d(X_u,x)\leq \varphi(u)$ for some $u\in (n_k,n_{k+1}]$}\right)\\
&\leq P_x\left(\text{$d(X_u,x)\leq \left(c/c_3\right)^{1/d_3}\varphi(n_k)$ for some $u\in (n_k,n_{k+1}]$}\right).
\end{split}
\end{equation}

We first give an upper bound of the last expression above by using Lemma \ref{prop:hit-prob}.
By (\ref{eq:scale-growth-inverse}), 
\begin{equation}\label{eq:comp-2}
\frac{\phi^{-1}((c-1)^2u)}{\phi^{-1}(u)}
\geq c_3^{1/d_3}\left(\frac{(c-1)^2u}{u}\right)^{1/d_3}=c_3^{1/d_3}(c-1)^{2/d_3}
\end{equation}
for any $u>0$. Since $g(t)\rightarrow 0$ as $t\rightarrow\infty$, 
there exists  $T_c>0$ such that 
$$g(u)\leq \left\{\frac{c_3^2(c-1)^2}{c}\right\}^{1/d_3} \quad \text{for all $u\geq T_c$.}$$
By this inequality and  (\ref{eq:comp-2}), we obtain
\begin{equation}\label{eq:comp-1}
\left(\frac{c}{c_3}\right)^{1/d_3}\varphi(u)
=\left(\frac{c}{c_3}\right)^{1/d_3}\phi^{-1}(u)g(u)
\leq \phi^{-1}((c-1)^2u) \quad \text{for all $u\geq T_c$}.
\end{equation}

We now suppose that   $t\geq T_c$. 
Since $n_k\geq t\geq T_c$ and 
$$(c-1)^2n_k=(c-1)(n_{k+1}-n_k)\leq n_k,$$ 
we see from (\ref{eq:comp-1}) that 
$$\left(\frac{c}{c_3}\right)^{1/d_3}\varphi(n_k)\leq \phi^{-1}((c-1)^2n_k)
=\phi^{-1}((c-1)(n_{k+1}-n_k))\leq \phi^{-1}(n_k).$$
Hence by Lemma \ref{prop:hit-prob}, we get
\begin{equation}\label{eq:hit-prob-est}
\begin{split}
&P_x\left(\text{$d(X_u,x)\leq \left(c/c_3\right)^{1/d_3}\varphi(n_k)$ for some $u\in (n_k,n_{k+1}]$}\right)\\
&\leq K_1\frac{V(x,\left(c/c_3\right)^{1/d_3}\varphi(n_k))}{\phi(\left(c/c_3\right)^{1/d_3}\varphi(n_k))}
\int_{n_k}^{n_{k+1}+(c-1)(n_{k+1}-n_k)}\frac{{\rm d}u}{V(x,\phi^{-1}(u))}.
\end{split}
\end{equation}

We next evaluate the last expression above. Since
$$
\int_{n_{k+1}}^{n_{k+1}+(c-1)(n_{k+1}-n_k)}\frac{{\rm d}u}{V(x,\phi^{-1}(u))}
\leq 
\frac{(c-1)(n_{k+1}-n_k)}{V(x,\phi^{-1}(n_{k+1}))}
\leq (c-1)\int_{n_k}^{n_{k+1}}\frac{{\rm d}u}{V(x,\phi^{-1}(u))},
$$
we obtain 
\begin{equation}
\begin{split}\label{eq:int-comp-0}
&\int_{n_k}^{n_{k+1}+(c-1)(n_{k+1}-n_k)}\frac{{\rm d}u}{V(x,\phi^{-1}(u))}\\
&=\int_{n_k}^{n_{k+1}}\frac{{\rm d}u}{V(x,\phi^{-1}(u))}
+\int_{n_{k+1}}^{n_{k+1}+(c-1)(n_{k+1}-n_k)}\frac{{\rm d}u}{V(x,\phi^{-1}(u))}
\leq c\int_{n_k}^{n_{k+1}}\frac{{\rm d}u}{V(x,\phi^{-1}(u))}.
\end{split}
\end{equation}
If $n_k<u\leq n_{k+1}$, then $\phi^{-1}(n_k)\leq \phi^{-1}(u)$ and $g(n_k)\leq g(u)/R_{c,t}$ so that 
$$\varphi(n_k)=\phi^{-1}(n_k)g(n_k)\leq \frac{1}{R_{c,t}}\phi^{-1}(u)g(u)=\frac{1}{R_{c,t}}\varphi(u).$$
Therefore, we have
$$V\left(x,\left(\frac{c}{c_3}\right)^{1/d_3}\varphi(n_k)\right)
\leq V\left(x,\left(\frac{c}{c_3}\right)^{1/d_3}\frac{1}{R_{c,t}}\varphi(u)\right)
\leq c_2\left(\frac{c}{c_3}\right)^{d_2/d_3}\frac{1}{(R_{c,t})^{d_2}}V\left(x,\varphi(u)\right)$$
by (\ref{eq:volume-growth}). 
Since  $\phi(\left(c/c_3\right)^{1/d_3}\varphi(n_k))\geq \phi(\varphi(u))$ by (\ref{eq:comp-radius}), 
we see from (\ref{eq:int-comp-0}) that 
\begin{equation*}
\begin{split}
&\frac{V(x,\left(c/c_3\right)^{1/d_3}\varphi(n_k))}{\phi(\left(c/c_3\right)^{1/d_3}\varphi(n_k))}
\int_{n_k}^{n_{k+1}+(c-1)(n_{k+1}-n_k)}\frac{{\rm d}u}{V(x,\phi^{-1}(u))}\\
&\leq cc_2\left(\frac{c}{c_3}\right)^{d_2/d_3}\frac{1}{(R_{c,t})^{d_2}}
\int_{n_k}^{n_{k+1}}\frac{V(x,\varphi(u))}{\phi(\varphi(u))}\frac{{\rm d}u}{V(x,\phi^{-1}(u))}.
\end{split}
\end{equation*}
Combining this with (\ref{eq:hit-prob-0}) and  (\ref{eq:hit-prob-est}), 
we arrive at the inequality (\ref{eq:hit-prob-1}). 
\qed
\bigskip

We can finish the proof of (\ref{eq:sup}) by Lemma \ref{lem:hit-prob-1}.  
In fact, it follows by   (\ref{eq:sup-est-1}) and (\ref{eq:hit-prob-1}) that 
\begin{equation*}
\begin{split}
q_{\varphi}(t,x)
&\leq \frac{K_1c_2}{c_3^{d_2/d_3}} \frac{c^{1+d_2/d_3}}{(R_{c,t})^{d_2}}
\sum_{k=0}^{\infty}\int_{n_k}^{n_{k+1}}\frac{V(x,\varphi(u))}{\phi(\varphi(u))}\frac{{\rm d}u}{V(x,\phi^{-1}(u))}\\
&=\frac{K_1c_2}{c_3^{d_2/d_3}} \frac{c^{1+d_2/d_3}}{(R_{c,t})^{d_2}}
\int_t^{\infty}\frac{V(x,\varphi(u))}{\phi(\varphi(u))}\frac{{\rm d}u}{V(x,\phi^{-1}(u))},
\end{split}
\end{equation*}
and thus 
$$\frac{q_{\varphi}(t,x)}{\displaystyle \int_t^{\infty}\frac{V(x,\varphi(u))}{\phi(\varphi(u))}\frac{{\rm d}u}{V(x,\phi^{-1}(u))}}
\leq \frac{K_1c_2}{c_3^{d_2/d_3}} \frac{c^{1+d_2/d_3}}{(R_{c,t})^{d_2}}.$$
Since $\lim_{c\rightarrow 1+0}(\lim_{t\rightarrow \infty}R_{c,t})=1$ by assumption, 
we get (\ref{eq:sup}) by letting $t\rightarrow \infty$ and then $c\rightarrow 1+0$ in the inequality above.
\bigskip

\subsection{Proof of (\ref{eq:inf})}\label{subsect:proof-2}
Fix positive constants $t$, $k$, and $l$ with $1<l<k<2$ 
and define a sequence $\{n_m\}_{m=0}^{\infty}$ by 
$$n_0=t, \ n_{2m+1}=kn_{2m}, \ n_{2m+2}=ln_{2m+1} \ (m\geq 0).$$
If $n_{2m}\leq u\leq n_{2m+1}$, then a calculation similar to (\ref{eq:comp-radius}) shows that  
\begin{equation}\label{eq:comp-radius-1}
\varphi(u)\geq \left(\frac{c_3}{k}\right)^{1/d_3}\varphi(n_{2m+1}).
\end{equation}
We now define the event $A_{2m} \ (m=0,1,2,\dots)$ by
\begin{equation}\label{eq:event}
A_{2m}=\left\{\text{$d(X_u,x)\leq \left(c_3/k\right)^{1/d_3}\varphi(n_{2m+1})$ for some $u\in (n_{2m},n_{2m+1}]$}\right\}.
\end{equation}
Then (\ref{eq:comp-radius-1}) yields that 
$$A_{2m}\subset \left\{\text{$d(X_u,x)\leq \varphi(u)$ for some $u\in (n_{2m},n_{2m+1}]$}\right\}$$ 
and hence
\begin{equation}\label{eq:inf-est-1}
\begin{split}
q_{\varphi}(t,x)
&=P_x\left(\bigcup_{m=0}^{\infty}
\left\{\text{$d(X_u,x)\leq \varphi(u)$ for some $u\in (n_m,n_{m+1}]$}\right\}\right)\\
&\geq P_x\left(\bigcup_{m=0}^{\infty}
\left\{\text{$d(X_u,x)\leq \varphi(u)$ for some $u\in (n_{2m},n_{2m+1}]$}\right\}\right)\\
&\geq P_x\left(\bigcup_{m=0}^{\infty}A_{2m}\right)
\geq \sum_{i=0}^{\infty}\left(P_x(A_{2i})-\sum_{j=i+1}^{\infty}P_x(A_{2i}\cap A_{2j})\right).
\end{split}
\end{equation}
The last inequality above is the so-called Bonferroni inequality 
(see, e.g., \cite[Exercise 1.6.10]{D10}).

Let $\kappa_t=\min\left\{1,g(t)^{d_3}/c_3\right\}$ for $t>0$. 
To calculate the last expression of (\ref{eq:inf-est-1}),
we first show

\begin{prop}\label{prop:upper-cap}
If $1<k<3/2$, $1<l<2-1/k$ and $\kappa_t<k(l-1)/2$, 
then there exists a positive constant $A(k,l)$ such that 
for any $i\geq 0$ and  $j\geq i+1$,
$$P_x(A_{2i}\cap A_{2j})\leq A(k,l)
P_x(A_{2i})\int_{n_{2j-1}}^{n_{2j+1}}\frac{V(x,\varphi(u))}{\phi(\varphi(u))}\frac{{\rm d}u}{V(x,\phi^{-1}(u))} \quad \text{for any $x\in M$.}$$
\end{prop}

The constant $A(k,l)$ will be given in (\ref{eq:const-a}) below.
For the proof of Proposition \ref{prop:upper-cap}, 
we calculate $P_x(A_{2i}\cap A_{2j})$ by using the strong Markov property. 
For $i\geq 0$, we define 
$$\sigma_{2i}=
\begin{cases}
\inf\left\{u\in (n_{2i}, n_{2i+1}] \mid d(X_u,X_0)\leq \left(c_3/k\right)^{1/d_3}\varphi(n_{2i+1})\right\}, & \text{if $\{\}\ne\emptyset$},\\
\infty, & \text{if $\{\}=\emptyset$}.
\end{cases}$$
If $j\geq i+1$, then by the strong Markov property,
\begin{equation}\label{eq:cap}
\begin{split}
&P_x(A_{2i}\cap A_{2j})\\
&=P_x\left(\sigma_{2i}\leq n_{2i+1}, \ 
\text{$d(X_u,x)\leq \left(c_3/k\right)^{1/d_3}\varphi(n_{2j+1})$ for some $u\in (n_{2j}, n_{2j+1}]$}\right)\\
&=E_x\left[F_j(X_{\sigma_{2i}}, n_{2j}-\sigma_{2i}, n_{2j+1}-\sigma_{2i}): 
\sigma_{2i}\leq n_{2i+1}\right]\\
&\leq P_x(\sigma_{2i}\leq n_{2i+1})\sup_{d(z,x)\leq \left(c_3/k\right)^{1/d_3}\varphi(n_{2i+1})}
F_j(z,n_{2j}-n_{2i+1},n_{2j+1}-n_{2i})\\
&=P_x(A_{2i})\sup_{d(z,x)\leq \left(c_3/k\right)^{1/d_3}\varphi(n_{2i+1})}
F_j(z,n_{2j}-n_{2i+1},n_{2j+1}-n_{2i}),
\end{split}
\end{equation}
where
\begin{equation}\label{eq:fj}
F_j(y,s_1,s_2)=P_y
\left(\text{$d(X_u,x)\leq \left(c_3/k\right)^{1/d_3}\varphi(n_{2j+1})$ for some $u\in (s_1, s_2]$}\right).
\end{equation}
To obtain an upper  bound of (\ref{eq:cap}), 
we use  the comparison of heat kernels:

\begin{lem}\label{lem:comp-heat}  {\rm (\cite[Lemma 4.2]{SW15}).} \ 
For any $t>0$ and $x,z\in M$ such that $d(x,z)\leq \phi^{-1}(t)$,  
\begin{equation}\label{eq:comp-heat-0}
p(t,z,y)\leq H_1p(t,x,y)
\end{equation}
for any $y\in M$, where $H_1=c_2c_42^{d_2+d_4}L_2/L_1$.
\end{lem}

\pf \ Suppose that $d(x,z)\leq \phi^{-1}(t)$. 
Since the heat kernel is symmetric by assumption, 
(\ref{estimate}) implies that 
$$p(t,z,y)=p(t,y,z)\leq \frac{L_2}{V(y,\phi^{-1}(t))}$$
for any $t>0$ and $y,z\in M$. 
 If $d(z,y)\leq \phi^{-1}(t)$, 
then we have $d(y,x)\leq d(y,z)+d(z,x)\leq 2\phi^{-1}(t)$ so that 
by (\ref{eq:volume-growth}) and \eqref{eq:scale-growth}, 
\begin{equation*}
\begin{split}
\frac{L_2}{V(y,\phi^{-1}(t))}
=\frac{L_2t}{V(y,\phi^{-1}(t))\phi(\phi^{-1}(t))}
&\leq \frac{L_2t}{V(y,d(y,x)/2))\phi(d(y,x)/2)}\\
&\leq c_2c_42^{d_2+d_4}L_2\frac{t}{V(y,d(y,x))\phi(d(y,x))}.
\end{split}
\end{equation*}
We thus get
\begin{equation*}
\begin{split}
p(t,z,y)
&\leq c_2c_42^{d_2+d_4}L_2\min\left\{\frac{1}{V(y,\phi^{-1}(t))}, \frac{t}{V(y,d(y,x))\phi(d(y,x))}\right\}\\
&\leq c_2c_42^{d_2+d_4}\frac{L_2}{L_1}p(t,y,x)
=H_1p(t,x,y).
\end{split}
\end{equation*}
The last equality above follows from the symmetry of the heat kernel. On the other hand, if $d(z,y)>\phi^{-1}(t)$, then $d(y,x)\leq 2d(y,z)$ so that 
we obtain (\ref{eq:comp-heat-0})  by  the same way as for the former case.
\qed
\bigskip

The following lemma gives an upper  bound of (\ref{eq:cap}).

\begin{lem}\label{lem:comp-start-1}
If $\kappa_t<k(l-1)$,  
then for any $x,z\in M$ with $d(x,z)\leq \left(c_3/k\right)^{1/d_3}\varphi(n_{2i+1})$,
\begin{equation}\label{eq:comp-start-1}
F_j(z,n_{2j}-n_{2i+1}, n_{2j+1}-n_{2i})\leq H_1F_j(x,n_{2j}-n_{2i+1}, n_{2j+1}-n_{2i}).
\end{equation}
\end{lem}

\pf \  
Suppose that  $\kappa_t<k(l-1)$. 
Then  for any $j\geq i+1$,
\begin{equation}\label{eq:comp-kappa}
\kappa_t n_{2i}\leq k(l-1)n_{2i}=n_{2i+2}-n_{2i+1}\leq n_{2j}-n_{2i+1}.
\end{equation}
On the other hand, since $g(u)\leq (c_3\kappa_t)^{1/d_3}$ for all $u\geq t$ and 
$$\phi^{-1}(u)\leq \phi^{-1}(\kappa_t u)\left(\frac{1}{c_3\kappa_t}\right)^{1/{d_3}}$$
by (\ref{eq:scale-growth-inverse}), we have 
\begin{equation}\label{eq:kappa}
\varphi(u)=\phi^{-1}(u)g(u)\leq \phi^{-1}(\kappa_tu)
\end{equation}
for any  $u\in [n_{2i},n_{2i+1}]$. 
Therefore, we obtain by (\ref{eq:comp-radius-1}) and  (\ref{eq:comp-kappa}), 
\begin{equation}\label{eq:dist-bound}
\left(\frac{c_3}{k}\right)^{1/d_3}\varphi(n_{2i+1})\leq \varphi(n_{2i})
\leq \phi^{-1}(\kappa_t n_{2i})\leq \phi^{-1}(n_{2j}-n_{2i+1}).
\end{equation}
This inequality and Lemma \ref{lem:comp-heat} show that if $d(x,z)\leq (c_3/k)^{1/d_3}\varphi(n_{2i+1})$, 
then 
\begin{equation}\label{eq:comp-heat}
p(n_{2j}-n_{2i+1},z,w)\leq H_1p(n_{2j}-n_{2i+1},x,w)
\end{equation}
for any $w\in M$.

By the Markov property, 
\begin{equation*}
\begin{split}
&F_j(z, n_{2j}-n_{2i+1}, n_{2j+1}-n_{2i})
=E_z\left[F_j({X_{n_{2j}-n_{2i+1}}}, 0,n_{2j+1}-n_{2i}-(n_{2j}-n_{2i+1}))\right]\\
&=\int_M p(n_{2j}-n_{2i+1},z,w)F_j(w, 0,n_{2j+1}-n_{2i}-(n_{2j}-n_{2i+1}))\,m({\rm d}w).
\end{split}
\end{equation*}
Then  (\ref{eq:comp-heat}) implies that for  any $x,z\in M$ with $d(x,z)\leq \left(c_3/k\right)^{1/d_3}\varphi(n_{2i+1})$, the last expression above is less than
\begin{equation*}
\begin{split}
&H_1
\int_M p(n_{2j}-n_{2i+1},x,w)F_j(w, 0,n_{2j+1}-n_{2i}-(n_{2j}-n_{2i+1}))\,m({\rm d}w)\\
&=H_1 F_j(x, n_{2j}-n_{2i+1}, n_{2j+1}-n_{2i}),
\end{split}
\end{equation*}
which is our assertion.
\qed
\bigskip

We further derive an  upper bound of (\ref{eq:comp-start-1}) 
by using the next lemma. 

\begin{lem}\label{lem:comp-prob-2}
Suppose that $1<k<3/2$ and $1<l<2-1/k$. If $\kappa_t<k(l-1)/2$,  then for any $x\in M$, $i\geq 0$ and $j\geq i+1$,
\begin{equation}\label{eq:comp-prob-2}
\begin{split}
&F_j(x,n_{2j}-n_{2i+1}, n_{2j+1}-n_{2i})\\
&\leq 
H_2\left(\frac{kl}{l-1}\right)^{d_2/d_3}
\frac{V(x,\left(c_3/k\right)^{1/d_3}\varphi(n_{2j+1}))}{\phi(\left(c_3/k\right)^{1/d_3}\varphi(n_{2j+1}))}
\int_{n_{2j-1}}^{n_{2j+1}}\frac{{\rm d}u}{V(x,\phi^{-1}(u))}
\end{split}
\end{equation}
for $H_2=3K_1c_2/c_3^{d_2/d_3}$.
\end{lem}

\pf \  We first check that Lemma \ref{prop:hit-prob} is applicable to the calculation of the left hand side of (\ref{eq:comp-prob-2}). Suppose that $\kappa_t<k(l-1)/2$.
Since
$$\kappa_tn_{2j}\leq \frac{1}{2}k(l-1)n_{2j}=\frac{1}{2}(n_{2j+2}-n_{2j}),$$
we get 
$$
\left(\frac{c_3}{k}\right)^{1/d_3}\varphi(n_{2j+1})\leq \phi^{-1}(\kappa_t n_{2j})\leq \phi^{-1}\left(\frac{1}{2}(n_{2j+2}-n_{2j+1})\right)
$$
by the same way as in (\ref{eq:dist-bound}). 
We further suppose that $1<k<3/2$ and $1<l< 2-1/k$. 
Then 
$$n_{2j+2}-n_{2j+1}\leq (n_{2j+1}-n_{2i})-(n_{2j}-n_{2i+1})$$
by direct calculation. 
Since $kl<2$ by assumption, we also have 
$$n_{2j+2}-n_{2j+1}=kl(n_{2j}-n_{2j-1})\leq 2(n_{2j}-n_{2i+1})$$
so that 
$$
\left(\frac{c_3}{k}\right)^{1/d_3}\varphi(n_{2j+1})\leq \phi^{-1}\left(\min\left\{n_{2j}-n_{2i+1}, \frac{1}{2}\{(n_{2j+1}-n_{2i})-(n_{2j}-n_{2i+1})\}\right\}\right).
$$
Accordingly, we can apply Lemma \ref{prop:hit-prob} 
to the left hand side of (\ref{eq:comp-prob-2}) to get 
\begin{equation}\label{eq:comp-prob-1}
\begin{split}
&F_j(x,n_{2j}-n_{2i+1}, n_{2j+1}-n_{2i})\\
&=P_x
\left(\text{$d(X_u,x)\leq \left(c_3/k\right)^{1/d_3}\varphi(n_{2j+1})$ for some $u\in (n_{2j}-n_{2i+1}, n_{2j+1}-n_{2i}]$}\right)\\
&\leq K_1
\frac{V(x,\left(c_3/k\right)^{1/d_3}\varphi(n_{2j+1}))}{\phi(\left(c_3/k\right)^{1/d_3}\varphi(n_{2j+1}))}
\int_{n_{2j}-n_{2i+1}}^{n_{2j+1}-n_{2i}+\left\{(n_{2j+1}-n_{2i})-(n_{2j}-n_{2i+1})\right\}/2}\frac{{\rm d}u}{V(x,\phi^{-1}(u))}. 
\end{split}
\end{equation}

We next evaluate the integral in the last expression above. 
For $j\geq i+1$, since 
$$n_{2j}-n_{2i+1}\geq n_{2j}-n_{2j-1}=(l-1)n_{2j-1}$$
and $n_{2j+1}=kln_{2j-1}$, 
it follows from (\ref{eq:volume-growth}) and (\ref{eq:scale-growth-inverse}) that 
$$
V(x,\phi^{-1}(n_{2j}-n_{2i+1}))
\geq V(x,\phi^{-1}((l-1)n_{2j-1}))\
\geq \frac{c_3^{d_2/d_3}}{c_2}\left(\frac{l-1}{kl}\right)^{d_2/d_3}V(x,\phi^{-1}(n_{2j+1})).
$$
Noting that 
$$n_{2j+1}-n_{2i}-(n_{2j}-n_{2i+1})\leq 2(n_{2j+1}-n_{2j-1}),$$
we obtain 
\begin{equation*}
\begin{split}
&\int_{n_{2j}-n_{2i+1}}^{n_{2j+1}-n_{2i}+\left\{n_{2j+1}-n_{2i}-(n_{2j}-n_{2i+1})\right\}/2}\frac{{\rm d}u}{V(x,\phi^{-1}(u))}
\leq \frac{3}{2}\frac{n_{2j+1}-n_{2i}-(n_{2j}-n_{2i+1})}{V(x,\phi^{-1}(n_{2j}-n_{2i+1}))}\\
&\leq \frac{3c_2}{c_3^{d_2/d_3}}\left(\frac{kl}{l-1}\right)^{d_2/d_3}\frac{n_{2j+1}-n_{2j-1}}{V(x,\phi^{-1}(n_{2j+1}))}
\leq \frac{3c_2}{c_3^{d_2/d_3}}\left(\frac{kl}{l-1}\right)^{d_2/d_3}\int_{n_{2j-1}}^{n_{2j+1}}\frac{{\rm d}u}{V(x,\phi^{-1}(u))}.
\end{split}
\end{equation*}
By this inequality with (\ref{eq:comp-prob-1}), we have (\ref{eq:comp-prob-2}). 
\qed
\bigskip

We see from Lemmas \ref{lem:comp-start-1} and \ref{lem:comp-prob-2} that if $1<k<3/2$, $1<l<2-1/k$ and $\kappa_t<k(l-1)/2$, 
then for any $x,z\in M$ with $d(x,z)\leq \left(c_3/k\right)^{1/d_3}\varphi(n_{2i+1})$,
\begin{equation}\label{eq:comp-start-3}
\begin{split}
&F_j(z,n_{2j}-n_{2i+1}, n_{2j+1}-n_{2i})\\
&\leq H_1 H_2\left(\frac{kl}{l-1}\right)^{d_2/d_3}
\frac{V(x,\left(c_3/k\right)^{1/d_3}\varphi(n_{2j+1}))}{\phi(\left(c_3/k\right)^{1/d_3}\varphi(n_{2j+1}))}
\int_{n_{2j-1}}^{n_{2j+1}}\frac{{\rm d}u}{V(x,\phi^{-1}(u))}.
\end{split}
\end{equation}

{\it Proof of Proposition {\rm \ref{prop:upper-cap}}.} \ 
Since $d_1>d_4$ by assumption, we have for any $x\in M$ and $0<s<r$,
$$\frac{V(x,s)}{\phi(s)}=\frac{V(x,s)}{V(x,r)}\frac{\phi(r)}{\phi(s)}\frac{V(x,r)}{\phi(r)}
\leq \frac{c_4}{c_1}\left(\frac{r}{s}\right)^{d_4-d_1}\frac{V(x,r)}{\phi(r)}
\leq \frac{c_4}{c_1}\frac{V(x,r)}{\phi(r)}$$
by (\ref{eq:volume-growth}) and (\ref{eq:scale-growth}). 
Hence if  $0<s<cr$ for some $c>1$, then 
\begin{equation}\label{eq:volume-scale}
\frac{V(x,s)}{\phi(s)}\leq \frac{c_4}{c_1}\frac{V(x,cr)}{\phi(cr)}
\leq \frac{c_2c_4}{c_1c_3}c^{d_2-d_3}\frac{V(x,r)}{\phi(r)}.
\end{equation}
For any  $u\in [n_{2j-1}, n_{2j+1}]$,  since we obtain
$$\left(\frac{c_3}{k}\right)^{1/d_3}\varphi(n_{2j+1})\leq l^{1/d_3}\varphi(u)$$
by the same way as in (\ref{eq:comp-radius-1}), 
it follows from (\ref{eq:volume-scale}) that 
$$\frac{V(x,\left(c_3/k\right)^{1/d_3}\varphi(n_{2j+1}))}{\phi(\left(c_3/k\right)^{1/d_3}\varphi(n_{2j+1}))}
\leq \frac{c_2c_4}{c_1c_3}l^{(d_2-d_3)/d_3}\frac{V(x,\varphi(u))}{\phi(\varphi(u))}.$$
Therefore, for any $t>0$ and  $k,l>0$ with $1<l<k<2$, 
and for any $j\geq 1$,
\begin{equation}\label{eq:int-comp}
\begin{split}
&\frac{V(x,\left(c_3/k\right)^{1/d_3}\varphi(n_{2j+1}))}{\phi(\left(c_3/k\right)^{1/d_3}\varphi(n_{2j+1}))}
\int_{n_{2j-1}}^{n_{2j+1}}\frac{{\rm d}u}{V(x,\phi^{-1}(u))}\\
&\leq 
\frac{c_2c_4}{c_1c_3}l^{(d_2-d_3)/d_3}
\int_{n_{2j-1}}^{n_{2j+1}}\frac{V(x,\varphi(u))}{\phi(\varphi(u))}\frac{{\rm d}u}{V(x,\phi^{-1}(u))}.
\end{split}
\end{equation}

We see by  (\ref{eq:comp-start-3}) and (\ref{eq:int-comp}) that  
if $1<k<3/2$, $1<l<2-1/k$ and $\kappa_t<k(l-1)/2$, then for any $x,z\in M$ with $d(x,z)\leq \left(c_3/k\right)^{1/d_3}\varphi(n_{2i+1})$,
\begin{equation}\label{eq:comp-start-4}
F_j(z,n_{2j}-n_{2i+1}, n_{2j+1}-n_{2i})
\leq A(k,l)
\int_{n_{2j-1}}^{n_{2j+1}}\frac{V(x,\varphi(u))}{\phi(\varphi(u))}\frac{{\rm d}u}{V(x,\phi^{-1}(u))}
\end{equation}
for
\begin{equation}\label{eq:const-a}
A(k,l)=H_1H_2\frac{c_2c_4}{c_1c_3}\left(\frac{kl}{l-1}\right)^{d_2/d_3}l^{(d_2-d_3)/d_3}.
\end{equation}
By (\ref{eq:cap}) and  (\ref{eq:comp-start-4}), 
the proof is complete. 
\qed 
\bigskip

Under the full conditions in Proposition \ref{prop:upper-cap}, we have  
\begin{equation*}
\begin{split}
\sum_{j=i+1}^{\infty}P_x(A_{2i}\cap A_{2j})
&\leq A(k,l)P_x(A_{2i})\int_{n_{2i+1}}^{\infty}\frac{V(x,\varphi(u))}{\phi(\varphi(u))}\frac{{\rm d}u}{V(x,\phi^{-1}(u))}\\
&\leq A(k,l)P_x(A_{2i})\int_t^{\infty}\frac{V(x,\varphi(u))}{\phi(\varphi(u))}\frac{{\rm d}u}{V(x,\phi^{-1}(u))}.
\end{split}
\end{equation*}
Therefore, (\ref{eq:inf-est-1}) shows that  
\begin{equation}\label{eq:inf-est-2}
q_{\varphi}(t,x)\geq \left(1-A(k,l)\int_t^{\infty}\frac{V(x,\varphi(u))}{\phi(\varphi(u))}
\frac{{\rm d}u}{V(x,\phi^{-1}(u))}\right)\sum_{i=0}^{\infty}P_x(A_{2i}).
\end{equation}

We next derive a lower bound of $P_x(A_{2i})$. 
\begin{prop}\label{prop:prob-lower}
If $g(t)<1$ and $\kappa_t<(k-1)/(2^{d_4}c_4)$, then for any $x\in M$ and $i\geq 0$, 
\begin{equation}\label{eq:prob-lower}
P_x(A_{2i})
\geq B(k,l)
\min\left\{R_{k,t}, \left(\frac{c_3}{l}\right)^{1/d_3}\right\}^{d_2-d_3}
\int_{n_{2i}}^{n_{2i+2}}\frac{V(x,\varphi(u))}{\phi(\varphi(u))}\frac{{\rm d}u}{V(x,\phi^{-1}(u))}
\end{equation}
for
$$B(k,l)=K_2\frac{c_1(c_3)^{2d_2/d_3}}{(c_2)^2c_4} \frac{k-1}{kl-1}\frac{1}{k^{2d_2/d_3-1}}.$$
\end{prop}

To show Proposition \ref{prop:prob-lower}, 
we first evaluate $P_x(A_{2i})$ by using Lemma \ref{lem:hit-lower}.  

\begin{lem}\label{lem:est-prob-l}
If $g(t)<1$ and $\kappa_t<(k-1)/(2^{d_4}c_4)$, 
then for any $x\in M$ and $i\geq 0$, 
\begin{equation}\label{eq:est-prob-l}
P_x(A_{2i})\geq K_2
\frac{V(x,(c_3/k)^{1/d_3}\varphi(n_{2i+1}))}{\phi((c_3/k)^{1/d_3}\varphi(n_{2i+1}))}
\int_{n_{2i}}^{n_{2i+1}}\frac{{\rm d}u}{V(x,\phi^{-1}(u))}.
\end{equation}
\end{lem}

\pf \ 
Let $A_{2i}$ be the event defined by (\ref{eq:event}). Then 
\begin{equation}\label{eq:prob-a}
P_x(A_{2i})=P_x\left(\text{$d(X_u,x)\leq \left(c_3/k\right)^{1/d_3}\varphi(n_{2i+1})$ 
for some $u\in (n_{2i},n_{2i+1}]$}\right).
\end{equation}
Suppose that  $g(t)<1$ and $\kappa_t<(k-1)/(2^{d_4}c_4)$. 
Then (\ref{eq:comp-radius-1}) shows that  
$$\phi((c_3/k)^{1/d_3}\varphi(n_{2i+1}))\leq \phi(\varphi(n_{2i}))\leq \phi(\phi^{-1}(n_{2i}))=n_{2i}.$$ 
We also have by (\ref{eq:comp-radius-1}) and (\ref{eq:kappa}),  
$$\phi(2(c_3/k)^{1/d_3}\varphi(n_{2i+1}))
\leq \phi(2\varphi(n_{2i}))\leq 
c_42^{d_4}\phi(\phi^{-1}(\kappa_t n_{2i}))
=c_42^{d_4}\kappa_tn_{2i}\leq n_{2i+1}-n_{2i}.$$
Hence we get (\ref{eq:est-prob-l}) by applying 
Lemma \ref{lem:hit-lower} to the right hand side of (\ref{eq:prob-a}).
\qed
\bigskip

In order to give a lower bound of the right hand side of (\ref{eq:est-prob-l}), we next show
\begin{lem}\label{lem:int-int}
For any $i\geq 0$,
\begin{equation}\label{eq:int-int}
\begin{split}
&\frac{V(x,(c_3/k)^{1/d_3}\varphi(n_{2i+1}))}{\phi((c_3/k)^{1/d_3}\varphi(n_{2i+1}))}
\int_{n_{2i}}^{n_{2i+1}}\frac{{\rm d}u}{V(x,\phi^{-1}(u))}\\
&\geq \frac{c_1(c_3)^{2d_2/d_3}}{(c_2)^2c_4} \frac{k-1}{kl-1}\frac{1}{k^{2d_2/d_3-1}}
\min\left\{R_{k,t}, \left(\frac{c_3}{l}\right)^{1/d_3}\right\}^{d_2-d_3}
\int_{n_{2i}}^{n_{2i+2}}\frac{V(x,\varphi(u))}{\phi(\varphi(u))}\frac{{\rm d}u}{V(x,\phi^{-1}(u))}.
\end{split}
\end{equation}
\end{lem}

\pf \ We begin by evaluating the integral in the left hand side of (\ref{eq:int-int}). 
By the definition of the sequence $\{n_m\}_{m=0}^{\infty}$,
\begin{equation}\label{eq:int-comp-1}
\begin{split}
\int_{n_{2i}}^{n_{2i+1}}\frac{{\rm d}u}{V(x,\phi^{-1}(u))}
&\geq \frac{n_{2i+1}-n_{2i}}{V(x,\phi^{-1}(n_{2i+1}))}
=\frac{k-1}{kl-1}\frac{n_{2i+2}-n_{2i}}{V(x,\phi^{-1}(n_{2i+1}))}\\
&=\frac{k-1}{kl-1}\left(\frac{n_{2i+2}-n_{2i+1}}{V(x,\phi^{-1}(n_{2i+1}))}+\frac{n_{2i+1}-n_{2i}}{V(x,\phi^{-1}(n_{2i+1}))}\right).
\end{split}
\end{equation}
Then 
$$\frac{n_{2i+2}-n_{2i+1}}{V(x,\phi^{-1}(n_{2i+1}))}\geq\int_{n_{2i+1}}^{n_{2i+2}}\frac{{\rm d}u}{V(x,\phi^{-1}(u))}$$
and 
$$\frac{n_{2i+1}-n_{2i}}{V(x,\phi^{-1}(n_{2i+1}))}\geq \frac{c_3^{d_2/d_3}}{c_2k^{d_2/d_3}}\frac{n_{2i+1}-n_{2i}}{V(x,\phi^{-1}(n_{2i}))}
\geq\frac{c_3^{d_2/d_3}}{c_2k^{d_2/d_3}}\int_{n_{2i}}^{n_{2i+1}}\frac{{\rm d}u}{V(x,\phi^{-1}(u))}$$
by (\ref{eq:volume-growth}) and (\ref{eq:scale-growth-inverse}). 
At the first  inequality above, we used the fact that  $c_3\leq 1$ and $k>1$. 
Since
\begin{equation*}
\begin{split}
&\frac{n_{2i+2}-n_{2i+1}}{V(x,\phi^{-1}(n_{2i+1}))}+\frac{n_{2i+1}-n_{2i}}{V(x,\phi^{-1}(n_{2i+1}))}\\
&\geq 
\int_{n_{2i+1}}^{n_{2i+2}}\frac{{\rm d}u}{V(x,\phi^{-1}(u))}
+\frac{c_3^{d_2/d_3}}{c_2k^{d_2/d_3}}\int_{n_{2i}}^{n_{2i+1}}\frac{{\rm d}u}{V(x,\phi^{-1}(u))}
\geq \frac{c_3^{d_2/d_3}}{c_2k^{d_2/d_3}}\int_{n_{2i}}^{n_{2i+2}}\frac{{\rm d}u}{V(x,\phi^{-1}(u))},
\end{split}
\end{equation*}
(\ref{eq:int-comp-1}) implies that 
\begin{equation}\label{eq:comp-int}
\int_{n_{2i}}^{n_{2i+1}}\frac{{\rm d}u}{V(x,\phi^{-1}(u))}
\geq  \frac{c_3^{d_2/d_3}}{c_2}\frac{k-1}{kl-1}\frac{1}{k^{d_2/d_3}}\int_{n_{2i}}^{n_{2i+2}}\frac{{\rm d}u}{V(x,\phi^{-1}(u))}.
\end{equation}

On account of (\ref{eq:comp-int}), the proof of (\ref{eq:int-int}) is completed by showing that 
\begin{equation}\label{eq:int-inside}
\begin{split}
&\frac{V(x,(c_3/k)^{1/d_3}\varphi(n_{2i+1}))}{\phi((c_3/k)^{1/d_3}\varphi(n_{2i+1}))}
\int_{n_{2i}}^{n_{2i+2}}\frac{{\rm d}u}{V(x,\phi^{-1}(u))}\\
&\geq  \frac{c_1c_3}{c_2c_4}\left(\frac{c_3}{k}\right)^{(d_2-d_3)/d_3}
\min\left\{R_{k,t}, \left(\frac{c_3}{l}\right)^{1/d_3}\right\}^{d_2-d_3}
\int_{n_{2i}}^{n_{2i+2}}\frac{V(x,\varphi(u))}{\phi(\varphi(u))}\frac{{\rm d}u}{V(x,\phi^{-1}(u))}.
\end{split}
\end{equation}
If $n_{2i}\leq u\leq n_{2i+1}$, then 
$$\varphi(n_{2i+1})=\phi^{-1}(n_{2i+1})g(n_{2i+1})
\geq R_{k,t}\phi^{-1}(u)g(u)=R_{k,t}\varphi(u).$$
On the other hand, if $n_{2i+1}\leq u\leq n_{2i+2}$, then 
$$\varphi(n_{2i+1})\geq \left(\frac{c_3}{l}\right)^{1/d_3}\varphi(u)$$
by (\ref{eq:comp-radius}). 
Hence for any $u\in [n_{2i},n_{2i+2}]$, we have 
$$\varphi(n_{2i+1})\geq \min\left\{R_{k,t}, \left(\frac{c_3}{l}\right)^{1/d_3}\right\}\varphi(u)$$
so that by (\ref{eq:volume-scale}), 
$$\frac{V(x,(c_3/k)^{1/d_3}\varphi(n_{2i+1}))}{\phi((c_3/k)^{1/d_3}\varphi(n_{2i+1}))}
\geq \frac{c_1c_3}{c_2c_4}\left(\frac{c_3}{k}\right)^{(d_2-d_3)/d_3}
\min\left\{R_{k,t}, \left(\frac{c_3}{l}\right)^{1/d_3}\right\}^{d_2-d_3}\frac{V(x,\varphi(u))}{\phi(\varphi(u))}.$$
This inequality  yields (\ref{eq:int-inside}).
\qed
\bigskip

{\it Proof of Proposition {\rm \ref{prop:prob-lower}}.} \ 
The assertion follows  by  Lemmas \ref{lem:est-prob-l} 
and \ref{lem:int-int}.   
\qed
\bigskip

We are now in a position to finish the proof of (\ref{eq:inf}). 
Under the full conditions in both Propositions \ref{prop:upper-cap} and \ref{prop:prob-lower}, 
since 
$$\sum_{i=0}^{\infty}P_x(A_{2i})\geq B(k,l)
\min\left\{R_{k,t}, \left(\frac{c_3}{l}\right)^{1/d_3}\right\}^{d_2-d_3}
\int_t^{\infty}\frac{V(x,\varphi(u))}{\phi(\varphi(u))}\frac{{\rm d}u}{V(x,\phi^{-1}(u))}$$
by (\ref{eq:prob-lower}),  we see from (\ref{eq:inf-est-2}) that   
\begin{equation*}
\begin{split}
q_{\varphi}(t,x)
& \geq B(k,l)\min\left\{R_{k,t}, \left(\frac{c_3}{l}\right)^{1/d_3}\right\}^{d_2-d_3}
\left(1-A(k,l)\int_t^{\infty}\frac{V(x,\varphi(u))}{\phi(\varphi(u))}\frac{{\rm d}u}{V(x,\phi^{-1}(u))}\right)
\\
&\times 
\int_t^{\infty}\frac{V(x,\varphi(u))}{\phi(\varphi(u))}\frac{{\rm d}u}{V(x,\phi^{-1}(u))}.
\end{split}
\end{equation*}
Namely,
\begin{equation*}
\begin{split}
&\frac{q_{\varphi}(t,x)}{\displaystyle \int_t^{\infty}\frac{V(x,\varphi(u))}{\phi(\varphi(u))}\frac{{\rm d}u}{V(x,\phi^{-1}(u))}}\\
&\geq B(k,l)\min\left\{R_{k,t}, \left(\frac{c_3}{l}\right)^{1/d_3}\right\}^{d_2-d_3}
\left(1-A(k,l)\int_t^{\infty}\frac{V(x,\varphi(u))}{\phi(\varphi(u))}\frac{{\rm d}u}{V(x,\phi^{-1}(u))}\right).
\end{split}
\end{equation*}
By letting first $t\rightarrow \infty$ and then $l\rightarrow 1+0$ and $k\rightarrow 1+0$, 
we arrive at (\ref{eq:inf}).

\appendix
\section{Proof of Theorem \ref{asymp-thm-critical}}\label{sect:appendix-proof}
In this appendix, we prove Theorem \ref{asymp-thm-critical}. 
More precisely, we show that 
\begin{equation}\label{eq:sup-critical}
\limsup_{t\rightarrow\infty}
\frac{q_{\varphi}(t,x)}
{\displaystyle \int_t^{\infty}\frac{{\rm d}s}{s|\log g(s)|}}
\leq \frac{L_2}{L_1}\frac{2^{d_2}c_2(c_{v,2})^2}{d_3(c_{v,1})^2}
\end{equation}
and 
\begin{equation}\label{eq:inf-critical}
\liminf_{t\rightarrow\infty}
\frac{q_{\varphi}(t,x)}
{\displaystyle \int_t^{\infty}\frac{{\rm d}s}{s|\log g(s)|}}
\geq \frac{L_1}{L_2}\frac{(c_{v,1})^2}{2^{d_4+1}d_4c_4(c_{v,2})^2}.
\end{equation}
Our approach here is similar to that of  Theorem  \ref{asymp-thm}. 
Throughout this appendix, we assume 
the full conditions in Theorem \ref{asymp-thm-critical}.
For simplicity, we also assume $N=\emptyset$ 
in Assumptions \ref{assum:heat-kernel} and \ref{assum:estimate}. 

\subsection{Hitting time distributions}
We first note that Lemmas \ref{lem:hit} and \ref{prop:hit-prob} are valid under the setting in this section. Using these lemmas, we have 
\begin{lem}\label{prop:hit-prob-critical}{\rm (\cite[Lemma 4.22]{SW15})}. \ 
Let $a$, $b$, $c$ and $r$ be positive constants. 
\begin{enumerate}
\item[{\rm (i)}] If $\phi(r)\leq a\wedge c$, 
then for any $x\in M$,
$$P_x\left(\text{$d(X_s,x)\leq r$ for some $s\in (a,b]$}\right)\leq K_3
\frac{\displaystyle \log\left(\frac{b+c}{a}\right)}{\displaystyle 1+\log\left(\frac{c}{\phi(r)}\right)},$$
where $K_3=K_1(c_{v,2}/c_{v,1})^2$.
\item[{\rm (ii)}] If $\phi(r)\leq a$ and $\phi(2r)\leq b-a$, 
then for any $x\in M$, 
\begin{equation*} 
\begin{split}
P_x\left(\text{$d(X_s,x)\leq r$ for some $s\in (a,b]$}\right)
\geq K_4
\frac{\displaystyle \log\left(\frac{b}{a}\right)}{\displaystyle 1+\log\left(\frac{b-a}{\phi(2r)}\right)},
\end{split}
\end{equation*}
where
$$K_4=\frac{L_1}{L_2}\frac{(c_{v,1})^2}{2^{d_4+1}c_4(c_{v,2})^2}.$$
\end{enumerate}
\end{lem}

\pf \ We first suppose  that $\phi(r)\leq a\wedge c$. 
Then by assumption and  (\ref{eq:hit-prob-est-1}),
$$\int_a^{b+c}P_x(d(X_u,x)\leq 2r)\,{\rm d}u
\leq \frac{2^{d_2}c_2c_{v,2}L_2}{c_{v,1}} \phi(r)\int_a^{b+c}\frac{{\rm d}u}{u}
=\frac{2^{d_2}c_2c_{v,2}L_2}{c_{v,1}}  \phi(r)\log\left(\frac{b+c}{a}\right).$$
We also see by assumption and   (\ref{eq:ball-inside}) that 
\begin{equation}\label{eq:crit-lower}
\begin{split}
P_y(d(X_u,y)\leq r)
\geq L_1\min\left\{1,\frac{V(y,r)}{V(y,\phi^{-1}(u))}\right\}
&\geq L_1\min\left\{1,\frac{c_{v,1}}{c_{v,2}}\frac{\phi(r)}{u}\right\}\\
&\geq \frac{c_{v,1}L_1}{c_{v,2}}\min\left\{1,\frac{\phi(r)}{u}\right\}
\end{split}
\end{equation}
for any $y\in M$ and $u>0$, and therefore
\begin{equation}\label{eq:inf-lower}
\begin{split}
\int_0^c\inf_{d(y,x)\leq r}P_y(d(X_u,y)\leq r)\,{\rm d}u
&\geq  \frac{c_{v,1}L_1}{c_{v,2}}\left(\int_0^{\phi(r)}\,{\rm d}u
+\phi(r)\int_{\phi(r)}^c \frac{{\rm d}u}{u}\right)\\
&=\frac{c_{v,1}L_1}{c_{v,2}}\phi(r)\left(1+\log\left(\frac{c}{\phi(r)}\right)\right).
\end{split}
\end{equation}
Hence the proof of (i) is completed by Lemma \ref{lem:hit}. 

We next suppose that $\phi(r)\leq a$ and $\phi(2r)\leq b-a$. 
Then by (\ref{eq:crit-lower}), 
$$\int_a^bP_x(d(X_u,x)\leq r)\,{\rm d}u
\geq \frac{c_{v,1}L_1}{c_{v,2}}\phi(r)\int_a^b\frac{{\rm d}u}{u}
=\frac{c_{v,1}L_1}{c_{v,2}}\phi(r)\log\left(\frac{b}{a}\right).$$
In a way similar to (\ref{eq:crit-lower}), we also have 
$$P_y(d(X_u,y)\leq 2r)
\leq \frac{c_{v,2}L_2}{c_{v,1}}\min\left\{1,\frac{\phi(2r)}{u}\right\}$$
for any $y\in M$ and $u>0$. 
Since $\phi(2r)\leq c_42^{d_4}\phi(r)$ by (\ref{eq:scale-growth}), we get 
\begin{equation*}
\begin{split}
\int_{0}^{b-a}\sup_{d(y,x)\leq r}P_y(d(X_u,y)\leq 2r)\,{\rm d}u
&\leq \frac{c_{v,2}c_42^{d_4}L_2}{c_{v,1}}\phi(r)\left(1+\log\left(\frac{b-a}{\phi(2r)}\right)\right)
\end{split}
\end{equation*}
by the same way as in \eqref{eq:inf-lower}. 
Hence the assertion (ii) follows by  Lemma \ref{lem:hit}. 
\qed

\subsection{Proof of (\ref{eq:sup-critical})} 
Since $g(t)\rightarrow 0$ as $t\rightarrow \infty$, 
we always take $t>0$ such that  $g(t)<1$. 
For a fixed constant $c\in (1,2)$, 
we define a sequence $\{n_k\}_{k=0}^{\infty}$ by $n_k=tc^k$ $(k\geq 0)$. 
In order to give an upper bound of the last expression in (\ref{eq:sup-est-1}), 
we show 

\begin{lem}\label{lem:hit-prob-est-critical}
For any $\varepsilon\in (0,d_3)$ and $c\in (1,2)$,  
there exists $T_{\varepsilon, c}>0$ such that 
for all $t\geq T_{\varepsilon, c}$, 
\begin{equation}\label{eq:hit-prob-est-critical}
\begin{split}
&P_x\left(\text{$d(X_u,x)\leq \varphi(u)$ for some $u\in (n_k,n_{k+1}]$}\right)\\
&\leq \frac{K_3}{d_3-\varepsilon}
\frac{c\log((c-1)^2+c)}{c-1}R_{c,t}\int_{n_k}^{n_{k+1}}\frac{{\rm d}u}{u|\log g(u)|}
\end{split}
\end{equation}
for any $x\in M$ and $k\geq 0$. 
\end{lem}

\pf \  We first note that 
\begin{equation}\label{eq:hit-prob-est-crit-0}
P_x\left(\text{$d(X_u,x)\leq \varphi(u)$ for some $u\in (n_k,n_{k+1}]$}\right)
\leq K_3
\frac{\displaystyle \log\left(\frac{n_{k+1}+(c-1)(n_{k+1}-n_k)}{n_k}\right)}
{\displaystyle 1+\log\left(\frac{(c-1)(n_{k+1}-n_k)}{\phi((c/c_3)^{1/d_3}\varphi(n_k))}\right)}.
\end{equation}
We can show this inequality by following the proof of Lemma \ref{lem:hit-prob-1} 
and using Lemma \ref{prop:hit-prob-critical} instead of Lemma \ref{prop:hit-prob}. 

We next evaluate the right hand side of (\ref{eq:hit-prob-est-crit-0}). 
Since $n_{k+1}=cn_k$, we have 
$$\log\left(\frac{n_{k+1}+(c-1)(n_{k+1}-n_k)}{n_k}\right)=\log((c-1)^2+c).$$
Moreover, since (\ref{eq:scale-growth}) implies that 
\begin{equation*}
\begin{split}
\phi\left(\left(\frac{c}{c_3}\right)^{1/d_3}\varphi(n_k)\right)
\leq c_4\left(\frac{c}{c_3}\right)^{d_4/d_3}\phi(\varphi(n_k))
&\leq c_4\left(\frac{c}{c_3}\right)^{d_4/d_3}\frac{1}{c_3}g(n_k)^{d_3}\phi(\phi^{-1}(n_k))\\
&=\frac{c^{d_4/d_3}c_4}{c_3^{1+d_4/d_3}}g(n_k)^{d_3}n_k,
\end{split}
\end{equation*}
we get 
\begin{equation*}
\begin{split}
\log\left(\frac{(c-1)(n_{k+1}-n_k)}{\phi((c/c_3)^{1/d_3}\varphi(n_k))}\right)
&\geq \log\left(\frac{c_3^{1+d_3/d_4}(c-1)^2 n_k}{c^{d_4/d_3} c_4g(n_k)^{d_3}n_k}\right)\\
&=\log\frac{(c-1)^2}{c^{d_4/d_3}}+\log \frac{c_3^{1+d_4/d_3}}{c_4}+d_3|\log g(n_k)|\\
&\geq \log\frac{(c-1)^2}{c^{d_4/d_3}}+\log \frac{c_3^{1+d_4/d_3}}{c_4}+\frac{d_3}{R_{c,t}}|\log g(n_{k+1})|.
\end{split}
\end{equation*}
In particular, for any $\varepsilon\in (0,d_3)$ and $c\in (1,2)$,
there exists $T_{\varepsilon, c}>0$ such that 
$$1+\log\frac{(c-1)^2}{c^{d_4/d_3}}+\log \frac{c_3^{1+d_4/d_3}}{c_4}
\geq -\frac{\varepsilon}{R_{c,t}}|\log g(n_{k+1})|
\quad \text{for all $t\geq T_{\varepsilon,c}$, }$$
and hence 
\begin{equation}\label{eq:hit-prob-est-critical-1}
\begin{split}
\frac{\displaystyle \log\left(\frac{n_{k+1}+(c-1)(n_{k+1}-n_k)}{n_k}\right)}
{\displaystyle 1+\log\left(\frac{(c-1)(n_{k+1}-n_k)}{\phi((c/c_3)^{1/d_3}\varphi(n_k))}\right)}
&\leq \frac{R_{c,t}}{d_3-\varepsilon}\frac{\log((c-1)^2+c)}{|\log g(n_{k+1})|}\\
&=\frac{R_{c,t}}{d_3-\varepsilon}
\frac{c\log((c-1)^2+c)}{c-1}\frac{1}{|\log g(n_{k+1})|}\frac{n_{k+1}-n_k}{n_{k+1}}\\
&\leq \frac{R_{c,t}}{d_3-\varepsilon}
\frac{c\log((c-1)^2+c)}{c-1}\int_{n_k}^{n_{k+1}}\frac{{\rm d}u}{u|\log g(u)|}
\end{split}
\end{equation}
for all $t\geq T_{\varepsilon,c}$. 
The proof is completed by  (\ref{eq:hit-prob-est-crit-0}) and (\ref{eq:hit-prob-est-critical-1}).
\qed
\bigskip

By (\ref{eq:sup-est-1}) and Lemma \ref{lem:hit-prob-est-critical}, 
we see that for any $\varepsilon\in (0,d_3)$ and $c\in (1,2)$, 
there exists $T_{\varepsilon, c}>0$ such that for all $t\geq T_{\varepsilon,c}$, 
\begin{equation}\label{eq:hit-prob-est-crit}
\begin{split}
q_{\varphi}(t,x)
\leq
\frac{K_3}{d_3-\varepsilon}\frac{c\log((c-1)^2+c)}{c-1}
R_{c,t}\int_t^{\infty}\frac{{\rm d}u}{u|\log g(u)|},
\end{split}
\end{equation}
and thus
$$\frac{q_{\varphi}(t,x)}{\displaystyle \int_t^{\infty}\frac{{\rm d}u}{u|\log g(u)|}}
\leq \frac{K_3}{d_3-\varepsilon}\frac{c\log((c-1)^2+c)}{c-1}R_{c,t}.$$
By letting $t\rightarrow\infty$ and then $
c\rightarrow 1+0$ and $\varepsilon\rightarrow +0$, 
we arrive at (\ref{eq:sup-critical}).

\subsection{Proof of  (\ref{eq:inf-critical})} 
Fix positive constants $t$, $k$ and $l$ with $1<l<k<2$. 
We define a sequence $\{n_m\}_{m=0}^{\infty}$ by 
$$n_0=t, \ n_{2m+1}=kn_{2m}, \ n_{2m+2}=ln_{2m+1} \ (m\geq 0).$$
Let $A_{2m}$ be the event defined by (\ref{eq:event}). 
By the same way as in the proof of (\ref{eq:inf}) 
(see Subsection \ref{subsect:proof-2}), 
we first give an upper bound of the probability $P_x(A_{2i}\cap A_{2j})$.

\begin{prop}\label{prop:upper-cap-1}
Suppose that $1<k<3/2$ and $1<l<2-1/k$. 
Then there exists $T_{k,l}>0$ such that for all $t\geq T_{k,l}$, 
$$P_x(A_{2i}\cap A_{2j})\leq A'(k,l)P_x(A_{2i})\int_{n_{2j-1}}^{n_{2j+1}}\frac{{\rm d}u}{u|\log g(u)|}$$
for any $x\in M$, $i,\geq 0$ and $j\geq i+1$, where
$$A'(k,l)=\frac{2H_1K_3}{d_3}\frac{kl}{kl-1}\log\left(\frac{3}{2}\left(\frac{kl}{l-1}\right)\right).$$
\end{prop}

Recall the notation $F_j(y,s_1,s_2)$ in \eqref{eq:fj}. 
For the proof of Proposition \ref{prop:upper-cap-1}, 
it is enough to show the next lemma.

\begin{lem}\label{lem:comp-prob-3}
Suppose that $1<k<3/2$ and $1<l<2-1/k$. 
Then there exists $T_{k,l}>0$ such that 
if $d(x,z)\leq (c_3/k)^{1/d_3}\varphi(n_{2i+1})$, 
then for all $t\geq T_{k,l}$, 
$$
F_j(z,n_{2j}-n_{2i+1},n_{2j+1}-n_{2i})
\leq A'(k,l)\int_{n_{2j-1}}^{n_{2j+1}}\frac{{\rm d}u}{u|\log g(u)|}
$$
for any $i\geq 0$ with $j\geq i+1$.
\end{lem}

\pf \ 
In a way similar to the proof of Lemma \ref{lem:comp-prob-2}, 
we can apply Lemmas \ref{prop:hit-prob-critical} and \ref{lem:comp-start-1}
to show that if $d(x,z)\leq (c_3/k)^{1/d_3}\varphi(n_{2i+1})$ and $\kappa_t<k(l-1)/2$, 
then for any $i\geq 0$ and  $j\geq i+1$, 
\begin{equation}\label{eq:comp-prob-1-critical}
\begin{split}
&F_j(z,n_{2j}-n_{2i+1},n_{2j+1}-n_{2i})
\leq H_1F_j(x,n_{2j}-n_{2i+1},n_{2j+1}-n_{2i})\\
&\leq H_1K_3
\frac{\displaystyle \log\left(\frac{n_{2j+1}-n_{2i}+(n_{2j+1}-n_{2j}+n_{2i+1}-n_{2i})/2}{n_{2j}-n_{2i+1}}\right)}
{\displaystyle 1+\log\left(\frac{n_{2j+1}-n_{2j}+n_{2i+1}-n_{2i}}{2\phi((c_3/k)^{1/d_3}\varphi(n_{2j+1}))}\right)}.
\end{split}
\end{equation}

Since
$$n_{2j+1}-n_{2i}\leq n_{2j+1}=kln_{2j-1}$$
and 
$$n_{2j}-n_{2i+1}\geq n_{2j}-n_{2j-1}=(l-1)n_{2j-1},$$
we have 
$$\frac{n_{2j+1}-n_{2i}}{n_{2j}-n_{2i+1}}\leq \frac{kl}{l-1}$$
so that 
$$\frac{n_{2j+1}-n_{2i}+(n_{2j+1}-n_{2j}+n_{2i+1}-n_{2i})/2}{n_{2j}-n_{2i+1}}
=\frac{3}{2}\left(\frac{n_{2j+1}-n_{2i}}{n_{2j}-n_{2i+1}}\right)-\frac{1}{2}
\leq \frac{3}{2}\left(\frac{kl}{l-1}\right).$$
On the other hand, we obtain 
\begin{equation*}
\begin{split}
\phi\left(\left(\frac{c_3}{k}\right)^{1/d_3}\varphi(n_{2j+1})\right)
&=\phi\left(\left(\frac{c_3}{k}\right)^{1/d_3}g(n_{2j+1})\phi^{-1}(n_{2j+1})\right)\\
&\leq \frac{1}{c_3}\left\{\left(\frac{c_3}{k}\right)^{1/d_3}g(n_{2j+1})\right\}^{d_3}n_{2j+1}
=g(n_{2j+1})^{d_3}n_{2j}
\end{split}
\end{equation*}
by  (\ref{eq:scale-growth}).
Noting that 
$$n_{2j+1}-n_{2j}+n_{2i+1}-n_{2i}\geq n_{2j+1}-n_{2j}=(k-1)n_{2j},$$
we get 
$$
\frac{n_{2j+1}-n_{2j}+n_{2i+1}-n_{2i}}{2\phi((c_3/k)^{1/d_3}\varphi(n_{2j+1}))}
\geq \frac{k-1}{2g(n_{2j+1})^{d_3}}.
$$
Hence if $g(t)<1$, then 
\begin{equation}\label{eq:comp-log}
\log\left(\frac{n_{2j+1}-n_{2j}+n_{2i+1}-n_{2i}}{2\phi((c_3/k)^{1/d_3}\varphi(n_{2j+1}))}\right)
\geq \log\left(\frac{k-1}{2g(n_{2j+1})^{d_3}}\right)=\log\frac{k-1}{2}+d_3|\log g(n_{2j+1})|.
\end{equation}

Here we note that  there exists $T_{k,l}>0$ such that for all $t\geq T_{k,l}$, 
we have $g(t)<1$, $\kappa_t<k(l-1)/2$ and 
$$1+\log\frac{k-1}{2}\geq -\frac{d_3}{2}|\log g(n_{2j+1})|.$$
We thus obtain for all $t\geq T_{k,l}$, 
$$1+\log\left(\frac{n_{2j+1}-n_{2j}+n_{2i+1}-n_{2i}}{2\phi((c_3/k)^{1/d_3}\varphi(n_{2j+1}))}\right)\geq \frac{d_3}{2} |\log g(n_{2j+1})|$$
by (\ref{eq:comp-log}), whence 
\begin{equation*}
\begin{split}
&\frac{\displaystyle \log\left(\frac{n_{2j+1}-n_{2i}+(n_{2j+1}-n_{2j}+n_{2i+1}-n_{2i})/2}{n_{2j}-n_{2i+1}}\right)}
{\displaystyle 1+\log\left(\frac{n_{2j+1}-n_{2j}+n_{2i+1}-n_{2i}}{2\phi((c_3/k)^{1/d_3}\varphi(n_{2j+1}))}\right)}
\leq 
\frac{2}{d_3}\frac{\displaystyle \log\left(\frac{3}{2}\left(\frac{kl}{l-1}\right)\right)}
{|\log g(n_{2j+1})|}\\
&=\frac{2}{d_3}\frac{kl}{kl-1}\log\left(\frac{3}{2}\left(\frac{kl}{l-1}\right)\right)\frac{1}
{|\log g(n_{2j+1})|}\frac{n_{2j+1}-n_{2j-1}}{n_{2j+1}}\\
&\leq \frac{2}{d_3}\frac{kl}{kl-1}\log\left(\frac{3}{2}\left(\frac{kl}{l-1}\right)\right)\int_{n_{2j-1}}^{n_{2j+1}}\frac{{\rm d}u}{u|\log g(u)|}.
\end{split}
\end{equation*}
By this inequality and  (\ref{eq:comp-prob-1-critical}), we complete the proof. 
\qed
\bigskip

{\it Proof of Proposition {\rm \ref{prop:upper-cap-1}}.} \ 
The assertion follows from (\ref{eq:cap}) and Lemma \ref{lem:comp-prob-3}. 
\qed
\bigskip
 
Under the full conditions in Proposition \ref{prop:upper-cap-1}, we have for all $t\geq T_{k,l}$,
$$
\sum_{j=i+1}^{\infty}P_x(A_{2i}\cap A_{2j})
\leq A'(k,l)P_x(A_{2i}) \int_{n_{2i+1}}^{\infty}\frac{{\rm d}u}{u|\log g(u)|}
\leq A'(k,l)P_x(A_{2i}) \int_t^{\infty}\frac{{\rm d}u}{u|\log g(u)|}
$$
so that by  (\ref{eq:inf-est-1}),
\begin{equation}\label{eq:inf-est-2-critical}
q_{\varphi}(t,x)\geq \left(1-A'(k,l)\int_t^{\infty}\frac{{\rm d}u}{u|\log g(u)|}\right)\sum_{i=0}^{\infty}P_x(A_{2i}).
\end{equation}

We next give a lower bound of $P_x(A_{2i})$. 
\begin{prop}\label{prop:prob-lower-critical}
For any $\varepsilon>0$, $k\in (1,2)$ and $l\in (1,2)$ with $k>l$, 
there exists $T_{\varepsilon,k,l}>0$ such that for all $t\geq T_{\varepsilon,k,l}$,  
\begin{equation*}
\begin{split}
P_x(A_{2i})\geq 
\frac{B'_{\varepsilon}(k,l)}{R_{k,t}} \int_{n_{2i}}^{n_{2i+2}}\frac{{\rm d}u}{u|\log g(u)|}
\end{split}
\end{equation*}
for any $x\in M$ and $i\geq 0$, where
$$B'_{\varepsilon}(k,l)=\frac{K_4}{d_4+\varepsilon}\frac{\log k}{kl-1} .$$
\end{prop}

\pf \ We assume that $g(t)<1$ and $\kappa_t<(k-1)/(2^{d_4}c_4)$. 
In a way similar to the proof of Lemma \ref{lem:est-prob-l}, 
we can apply Lemma \ref{prop:hit-prob-critical} to show that 
$$
P_x(A_{2i})\geq 
K_4
\frac{\displaystyle \log k}{\displaystyle 1+\log\left(\frac{(k-1)n_{2i}}{\phi(2(c_3/k)^{1/d_3}\varphi(n_{2i+1}))}\right)}.
$$
Since 
\begin{equation*}
\begin{split}
\phi\left(2\left(\frac{c_3}{k}\right)^{1/d_3}\varphi(n_{2i+1})\right)
&\geq 2^{d_3}c_3\phi\left(\left(\frac{c_3}{k}\right)^{1/d_3}\phi^{-1}(n_{2i+1})g(n_{2i+1})\right)\\
&\geq 2^{d_3}\frac{c_3^{1+d_4/d_3}}{c_4k^{d_4/d_3}}g(n_{2i+1})^{d_4}n_{2i+1}
\end{split}
\end{equation*} 
by \eqref{eq:scale-growth}, 
we have 
\begin{equation*}
\begin{split}
&\log\left(\frac{(k-1)n_{2i}}{\phi(2(c_3/k)^{1/d_3}\varphi(n_{2i+1}))}\right)
\leq \log\left(\frac{c_4k^{d_4/d_3}(k-1)n_{2i}}{2^{d_3}c_3^{1+d_4/d_3}g(n_{2i+1})^{d_4}n_{2i+1}}\right)\\
&=\log\left(\frac{c_4}{2^{d_3}c_3^{1+d_4/d_3}}\right)+\log k^{d_4/d_3-1}(k-1)+d_4|\log g(n_{2i+1})|.
\end{split}
\end{equation*}
Note that  for any $\varepsilon>0$, $k\in (1,2)$ and $l\in (1,2)$ with $k>l$, 
there exists $T_{\varepsilon, k,l}>0$ such that for all $t\geq T_{\varepsilon, k,l}$, 
we obtain $g(t)<1$,  $\kappa_t<(k-1)/(2^{d_4}c_4)$ and 
$$1+\log\left(\frac{c_4}{2^{d_3}c_3^{1+d_4/d_3}}\right)+\log k^{d_4/d_3-1}(k-1)
\leq \varepsilon|\log g(n_{2i+1})|.$$
Hence for all  $t\geq T_{\varepsilon,k,l}$,  
$$1+\log\left(\frac{c_4}{2^{d_3}c_3^{1+d_4/d_3}}\right)+\log k^{d_4/d_3-1}(k-1)+d_4|\log g(n_{2i+1})|
\leq (d_4+\varepsilon)|\log g(n_{2i+1})|$$
so that 
\begin{equation}\label{eq:est-prob-lower-critical}
P_x(A_{2i})
\geq \frac{K_4}{d_4+\varepsilon}\frac{\log k}{|\log g(n_{2i+1})|}
=B'_{\varepsilon}(k,l)
\frac{1}{|\log g(n_{2i+1})|}\left(\frac{n_{2i+1}-n_{2i}}{n_{2i}}+\frac{n_{2i+2}-n_{2i+1}}{n_{2i}}\right).
\end{equation}

In what follows, we assume that $t\geq T_{\varepsilon,k,l}$. 
Let us evaluate the last expression of (\ref{eq:est-prob-lower-critical}).
By the definition of $R_{k,t}$,  
\begin{equation*}
\begin{split}
\frac{1}{|\log g(n_{2i+1})|}\frac{n_{2i+1}-n_{2i}}{n_{2i}}
&=\frac{|\log g(n_{2i})|}{|\log g(n_{2i+1})|}
\frac{n_{2i+1}-n_{2i}}{n_{2i}|\log g(n_{2i+1})|}\\
&\geq \frac{1}{R_{k,t}}\int_{n_{2i}}^{n_{2i+1}}\frac{{\rm d}u}{u|\log g(u)|}.
\end{split}
\end{equation*}
Since
$$\frac{1}{|\log g(n_{2i+1})|}\frac{n_{2i+2}-n_{2i+1}}{n_{2i}}
\geq \int_{n_{2i+1}}^{n_{2i+2}}\frac{{\rm d}u}{u|\log g(u)|},$$
we get 
$$\frac{1}{|\log g(n_{2i+1})|}\left(\frac{n_{2i+1}-n_{2i}}{n_{2i}}+\frac{n_{2i+2}-n_{2i+1}}{n_{2i}}\right)
\geq \frac{1}{R_{k,t}}\int_{n_{2i}}^{n_{2i+2}}\frac{{\rm d}u}{u|\log g(u)|}.$$
Combining this with (\ref{eq:est-prob-lower-critical}), we complete the proof. 
\qed
\bigskip
 
We are now in a position to finish the proof of (\ref{eq:inf-critical}). 
Under the full conditions in Propositions \ref{prop:upper-cap-1} and \ref{prop:prob-lower-critical}, 
we obtain for all $t\geq \max\{T_{k,l}, T_{\varepsilon,k,l}\}$,
$$\sum_{i=0}^{\infty}P_x(A_{2i})\geq \frac{B'_{\varepsilon}(k,l)}{R_{k,t}}
\int_t^{\infty}\frac{{\rm d}u}{u|\log g(u)|}.$$
Therefore, it follows by (\ref{eq:inf-est-2-critical}) that 
\begin{equation*}
\begin{split}
\frac{q_{\varphi}(t,x)}{\displaystyle \int_t^{\infty}\frac{{\rm d}u}{u|\log g(u)|}}
\geq \frac{B'_{\varepsilon}(k,l)}{R_{k,t}}\left(1-A'(k,l)\int_t^{\infty}\frac{{\rm d}u}{u|\log g(u)|}\right).
\end{split}
\end{equation*}
By letting first $t\rightarrow \infty$ and then $l\rightarrow 1+0$, 
$k\rightarrow 1+0$, and $\varepsilon\rightarrow+0$, 
we get (\ref{eq:inf-critical}).

\section{Derivation of (\ref{eq:dirichlet-sub}) and  (\ref{eq:heat-sub})}\label{sect:sub}
In this appendix, we show (\ref{eq:dirichlet-sub}) and  (\ref{eq:heat-sub}) in Example \ref{ex:sub-diff} above. 

\subsection{Subordination}
We first recall the notion of subordinators  according to \cite{B96} and \cite{Sa13}. 
An increasing L\'evy process on  $[0,\infty)$ is called a {\it subordinator}. 
By \cite[Theorem 21.5]{Sa13}, subordinators are characterized by the following Laplace transform: 
if we denote by $\pi_t({\rm d}s)$ the transition function of a subordinator, 
then 
$$\int_0^{\infty}e^{-\lambda s}\,\pi_t({\rm d}s)=e^{-t\psi(\lambda)} \quad (\lambda>0, \ t>0)$$
for 
$$\psi(\lambda)=b\lambda+\int_{(0,\infty)}(1-e^{-\lambda s})\nu({\rm d}s).$$
Here $b$ is a nonnegative constant and $\nu$ is a positive Radon measure on $(0,\infty)$ 
such that 
$$\int_{(0,\infty)}(s\wedge 1)\nu({\rm d}s)<\infty.$$
For $\gamma\in (0,1)$, we say that a subordinator $\{\tau_t\}_{t\geq 0}$ is  {\it $\gamma$-stable} 
if $\psi(\lambda)=\lambda^{\gamma}$, that is, 
\begin{equation}\label{eq:sub-stable}
b=0, \quad \nu({\rm d}s)=\frac{\gamma}{\Gamma(1-\gamma)}\frac{{\rm d}s}{s^{1+\gamma}}
\end{equation}
(see, e.g., \cite[Examples 21.7 and 24.12]{Sa13}).

We next introduce the subordination of symmetric Markov processes. 
Let $(M,d)$ be a locally compact separable metric space and $m$ 
a positive Radon measure on $M$ with full support.  
Let ${\bf M}=(\Omega, {\cal F}, \{X_t\}_{t\geq 0}, \{P_x\}_{x\in M})$ 
be an $m$-symmetric Hunt process on $M$ such that 
the corresponding Dirichlet form $({\cal E}, {\cal F})$ is regular on $L^2(M;m)$.
Fix a subordinator $\{\tau_t\}_{t\geq 0}$ defined on $(\Omega, {\cal F})$ 
such that it is independent of $\{X_t\}_{t\geq 0}$ 
under $P_x$ for every $x\in M$.  
Let ${\bf M}^{(1)}=(\{Y_t\}_{t\geq 0}, \{P_x\}_{x\in M})$ be a subordinated process of ${\bf M}$ 
defined by 
$$Y_t=X_{\tau_t} \quad \text{for $t\geq 0$}.$$
If we denote by $p(t,x,{\rm d}y)$ the transition function of ${\bf M}$, 
then the transition function of ${\bf M}^{(1)}$ is given by 
\begin{equation}\label{eq:sub-heat}
q(t,x,{\rm d}y)=\int_{0}^{\infty}p(s,x, {\rm d}y)\,\pi_t({\rm d}s)
\end{equation}
(\cite[Theorem 30.1]{Sa13}).
According to \cite[Theorem 2.1]{O02}, 
${\bf M}^{(1)}$ associates a regular Dirichlet form $({\cal E}^{(1)}, {\cal F}^{(1)})$ on $L^2(M;m)$
such that if $b>0$, then ${\cal F}^{(1)}={\cal F}$ and 
\begin{equation*}
{\cal E}^{(1)}(u,v)
=b{\cal E}(u,v)+\iint_{M\times M\setminus{\rm diag}}(u(x)-u(y))(v(x)-v(y))\,J({\rm d}x{\rm d}y)+\int_M u(x)v(x)\,k({\rm d}x)
\end{equation*}
for 
\begin{equation}\label{eq:j-kernel}
J({\rm d}x{\rm d}y)=\frac{1}{2}{\bf 1}_{\{x\ne y\}}m({\rm d}x)\int_{(0,\infty)}\nu({\rm d}s)p(s,x,{\rm d}y)
\end{equation}
and 
$$\ k({\rm d}x)=m({\rm d}x)\int_{(0,\infty)}(1-p(s,x,M))\,\nu({\rm d}s).$$
On the other hand, if  $b=0$, then the form ${\cal E}^{(1)}$ is the same as before 
and 
$$
{\cal F}^{(1)}=
\left\{u\in L^2(M;m) \mid \int_{(0,\infty)}\left(\int_M (u(x)-T_su(x))u(x)\,m({\rm d}x)\right)\nu({\rm d}s)<\infty\right\}.
$$
We note that ${\cal F}\subset {\cal F}^{(1)}$.

\subsection{Derivation of (\ref{eq:dirichlet-sub})}
\label{subsect:j-kernel}
Throughout this subsection, 
we assume the full conditions in Example \ref{ex:sub-diff}. 
We take the measure $\nu$ as in (\ref{eq:sub-stable}) for some $\gamma\in (0,1)$. 
We recall that  
$$\phi(r)=r^{\gamma\beta_1}{\bf 1}_{\{r<1\}}+r^{\gamma\beta_2}{\bf 1}_{\{r\geq 1\}}.$$
On account of (\ref{eq:j-kernel}), it is enough for the proof of (\ref{eq:dirichlet-sub}) to show that  
\begin{equation}\label{eq:jump-kernel}
\int_0^{\infty}p(s,x,y)\,\nu({\rm d}s)
\asymp \frac{1}{V(x,d(x,y))\phi(d(x,y))}.
\end{equation}

We only show the upper bound of the left side in (\ref{eq:jump-kernel}) because 
the lower bound follows by the same way. 
Define 
$$({\rm A})=\int_0^{1\vee d(x,y)}p(s,x,y)\,\nu({\rm d}s) \quad \text{and} \quad 
({\rm B})=\int_{1\vee d(x,y)}^{\infty}p(s,x,y)\,\nu({\rm d}s).$$
Then by the change of variables $u=d(x,y)^{\beta_1}/s$,
\begin{equation*}
\begin{split}
({\rm A})
&\lesssim 
\int_0^{1\vee d(x,y)}\frac{1}{V(x,s^{1/\beta_1})}\exp\left\{-C_2\left(\frac{d(x,y)^{\beta_1}}{s}\right)^{1/(\beta_1-1)}\right\}
\frac{{\rm d}s}{s^{1+\gamma}}\\
&\asymp 
\frac{1}{d(x,y)^{\gamma\beta_1}}\int_{d(x,y)^{\beta_1}\wedge d(x,y)^{{\beta_1-1}}}^{\infty}
\frac{1}{V(x,d(x,y)/u^{1/\beta_1})}e^{-C_2u^{1/(\beta_1-1)}}u^{\gamma-1}\,{\rm d}u=:({\rm A}').
\end{split}
\end{equation*}
Since 
\begin{equation}\label{eq:v-comp}
\begin{split}
\frac{V(x,d(x,y)/{u^{1/\beta_1}})}{V(x,d(x,y))}
\geq 
\begin{cases}
\displaystyle \frac{c_3}{u^{d_1/\beta_1}}, & (0<u<1),\\
\displaystyle \frac{1}{c_4u^{d_2/\beta_1}}, & (u\geq 1)
\end{cases}
\end{split}
\end{equation}
by (\ref{eq:volume-growth}), we have
\begin{equation*}
\begin{split}
({\rm A}')
&\lesssim \frac{1}{V(x,d(x,y))d(x,y)^{\gamma\beta_1}}\\
&\times \int_{d(x,y)^{\beta_1}\wedge d(x,y)^{{\beta_1-1}}}^{\infty}
\left(u^{d_1/\beta_1}{\bf 1}_{\{u<1\}}+u^{d_2/\beta_1}{\bf 1}_{\{u\geq 1\}}\right)
\exp\left(-C_2u^{1/(\beta_1-1)}\right)u^{\gamma-1}\,{\rm d}u.
\end{split}
\end{equation*}
By the same way, we obtain \begin{equation*}
\begin{split}
({\rm B})
&\lesssim \frac{1}{V(x,d(x,y))d(x,y)^{\gamma\beta_2}}\\
&\times \int_0^{d(x,y)^{\beta_2}\wedge d(x,y)^{{\beta_2-1}}}
\left(u^{d_1/\beta_2}{\bf 1}_{\{u<1\}}+u^{d_2/\beta_2}{\bf 1}_{\{u\geq 1\}}\right)
\exp\left(-C_4u^{1/(\beta_2-1)}\right)u^{\gamma-1}\,{\rm d}u
\end{split}
\end{equation*}
so that 
$$\int_0^{\infty}p(s,x,y)\,\nu({\rm d}s)=({\rm A})+({\rm B})\lesssim \frac{1}{V(x,d(x,y))\phi(d(x,y))}.$$

\subsection{Derivation of (\ref{eq:heat-sub})}
Throughout this subsection, 
we keep the same setting as in Subsection \ref{subsect:j-kernel}. 
Let $\pi_t(s)$ be the density of the transition function for the $\gamma$-stable subordinator. 
Then the following relations hold (see, e.g., \cite[Theorem 3.1]{BSS03} and \cite[Remark 14.18]{Sa13}):
\begin{itemize}
\item For each $t>0$, $\pi_t(s)$ is a bounded continuous function on $(0,\infty)$ 
such that 
\begin{equation}\label{eq:scale}
\pi_t(s)=\frac{1}{t^{1/\gamma}}\pi_1\left(\frac{s}{t^{1/\gamma}}\right)
\end{equation}
for any $s>0$ and $t>0$;
\item There exists $c>0$ such that 
\begin{equation}\label{eq:upper-pi}
\pi_t(s)\leq c\frac{t}{s^{1+\gamma}}\exp\left(-\frac{t}{s^{\gamma}}\right)\leq c\frac{t}{s^{1+\gamma}}
\end{equation}
for any $s>0$ and $t>0$;
\item There exists $c>0$ such that if $s\geq t^{1/\gamma}$, then 
\begin{equation}\label{eq:lower-pi}
\pi_t(s)\geq c\frac{t}{s^{1+\gamma}}.
\end{equation}
\end{itemize}
By (\ref{eq:sub-heat}), the $\gamma$-stable subordinated diffusion process 
admits the  heat kernel $q(t,x,y)$ such that 
$$q(t,x,y)=\int_0^{\infty}p(s,x,y)\pi_t(s)\,{\rm d}s.$$

To get the upper bound of \eqref{eq:heat-sub}, 
we first show that 
\begin{equation}\label{eq:off-diag}
q(t,x,y)\lesssim \frac{t}{V(x,d(x,y))\phi(d(x,y))}.
\end{equation}
We divide $q(t,x,y)$ into 
$$({\rm I})=\int_0^{1\vee d(x,y)}p(s,x,y)\pi_t(s)\,{\rm d}s 
\quad \text{and}\quad 
({\rm II})=\int_{1\vee d(x,y)}^{\infty}p(s,x,y)\pi_t(s)\,{\rm d}s.$$
By (\ref{eq:upper-pi}) and the change of variables $s=d(x,y)^{\beta_1}/u$,
\begin{equation*}
\begin{split}
&({\rm I})
\lesssim t\int_0^{1\vee d(x,y)}\frac{1}{V(x,s^{1/\beta_1})}
\exp\left\{-C_2\left(\frac{d(x,y)^{\beta_1}}{s}\right)^{1/(\beta_1-1)}\right\}\frac{{\rm d}s}{s^{1+\gamma/2}}\\
&=\frac{t}{d(x,y)^{\gamma\beta_1/2}}\int_{d(x,y)^{\beta_1}\wedge d(x,y)^{\beta_1-1}}^{\infty}\frac{1}{V(x,d(x,y)/u^{1/\beta_1})}
\exp\left(-C_2u^{1/(\beta_1-1)}\right)u^{\gamma-1}\,{\rm d}u=:({\rm I})'.
\end{split}
\end{equation*}
Then by  (\ref{eq:v-comp}),
\begin{equation*}
\begin{split}
({\rm I})'
&\lesssim  
\frac{t}{V(x,d(x,y))d(x,y)^{\gamma\beta_1}}\\
&\times \int_{d(x,y)^{\beta_1}\wedge d(x,y)^{\beta_1-1}}^{\infty}
(u^{d_1/\beta_1}{\bf 1}_{\{u<1\}}+u^{d_2/\beta_1}{\bf 1}_{\{u\geq 1\}})
\exp\left(-C_2u^{1/(\beta_1-1)}\right)u^{\gamma-1}\,{\rm d}u.
\end{split}
\end{equation*}
By the same way, we obtain 
\begin{equation*}
\begin{split}
({\rm II})
&\lesssim  
\frac{t}{V(x,d(x,y))d(x,y)^{\gamma\beta_2}}\\
&\times \int_0^{d(x,y)^{\beta_2}\wedge d(x,y)^{\beta_2-1}}
(u^{d_1/\beta_2}{\bf 1}_{\{u<1\}}+u^{d_2/\beta_2}{\bf 1}_{\{u\geq 1\}})
\exp\left(-C_4u^{1/(\beta_2-1)}\right)u^{\gamma-1}\,{\rm d}u
\end{split}
\end{equation*}
so that \eqref{eq:off-diag} follows.

We next show that 
\begin{equation}\label{eq:on-diag}
q(t,x,y)=({\rm I})+({\rm II})\lesssim \frac{1}{V(x,\phi^{-1}(t))}.
\end{equation}
Suppose first that $d(x,y)<1$. 
By \eqref{eq:upper-pi} and 
the change of variables $u=t/s^{\gamma}$,
$$
({\rm I})
\lesssim t\int_0^1\frac{1}{V(x,s^{1/\beta_1})}
\exp\left(-\frac{t}{s^{\gamma}}\right)\frac{{\rm d}s}{s^{1+\gamma}}
\asymp \int_t^{\infty}\frac{e^{-u}}{V(x,(t/u)^{1/(\gamma\beta_1)})}\,{\rm d}u=:({\rm I})'.
$$
Then by \eqref{eq:v-comp},
$$({\rm I})'\lesssim \frac{1}{V(x,t^{1/(\gamma\beta_1)})}
\int_t^{\infty}
(u^{d_1/(\gamma\beta_1)}{\bf 1}_{\{u<1\}}
+u^{d_2/(\gamma \beta_1)}{\bf 1}_{\{u\geq 1\}})e^{-u}\,{\rm d}u.$$
By the same way, we get 
\begin{equation*}
\begin{split}
({\rm II})
&\lesssim t\int_1^{\infty}\frac{1}{V(x,s^{1/\beta_2})}
\exp\left(-\frac{t}{s^{\gamma}}\right)\frac{{\rm d}s}{s^{1+\gamma}}
\asymp \int_0^t\frac{e^{-u}}{V(x,(t/u)^{1/(\gamma\beta_2)})}
\,{\rm d}u\\
&\lesssim \frac{1}{V(x,t^{1/(\gamma\beta_2)})}
\int_0^t
(u^{d_1/(\gamma\beta_2)}{\bf 1}_{\{u<1\}}
+u^{d_2/(\gamma \beta_2)}{\bf 1}_{\{u\geq 1\}})e^{-u}\,{\rm d}u.
\end{split}
\end{equation*}
Hence if $0<t<1$, then 
$$({\rm I})+({\rm II})
\lesssim t\int_0^1\frac{1}{V(x,s^{1/\beta_1})}
\exp\left(-\frac{t}{s^{\gamma}}\right)\frac{{\rm d}s}{s^{1+\gamma}}
\lesssim \frac{1}{V(x,t^{1/(\gamma\beta_1)})}.$$
On the other hand, if $t\geq 1$, then 
$$({\rm I})+({\rm II})
\lesssim \frac{1}{V(x,t^{1/(\gamma\beta_2)})}
\int_0^t
(u^{d_1/(\gamma\beta_2)}{\bf 1}_{\{u<1\}}
+u^{d_2/(\gamma \beta_2)}{\bf 1}_{\{u\geq 1\}})e^{-u}\,{\rm d}u
\asymp \frac{1}{V(x,t^{1/(\gamma\beta_2)})}.$$
Therefore,  \eqref{eq:on-diag} holds.

Suppose next that $d(x,y)\geq 1$. 
If $t<\phi(d(x,y))(=d(x,y)^{\gamma \beta_2})$, then \eqref{eq:on-diag} 
follows by (\ref{eq:off-diag}). 
In what follows, we assume that  $t\geq \phi(d(x,y))$.
By \eqref{eq:upper-pi} and the change of variables $u=t/s^{\gamma}$,
$$
({\rm I})
\lesssim t\int_0^{d(x,y)}\frac{1}{V(x,s^{1/\beta_1})}
\exp\left(-\frac{t}{s^{\gamma}}\right)\frac{{\rm d}s}{s^{1+\gamma}}
\asymp \int_{t/d(x,y)^{\gamma}}^{\infty}\frac{e^{-u}}{V(x,(t/u)^{1/(\gamma\beta_1)})}
\,{\rm d}u
=:({\rm I})''.
$$
Then by \eqref{eq:v-comp}, 
$$({\rm I})''\lesssim \frac{1}{V(x,t^{1/(\gamma\beta_1)})}\int_{t/d(x,y)^{\gamma}}^{\infty}
(u^{d_1/(\gamma\beta_1)}{\bf 1}_{\{u<1\}}+u^{d_2/(\gamma \beta_1)}{\bf 1}_{\{u\geq 1\}})e^{-u}\,{\rm d}u.$$
By the same way, we get 
$$({\rm II})\lesssim \frac{1}{V(x,t^{1/(\gamma\beta_2)})}\int_0^{t/d(x,y)^{\gamma}}
(u^{d_1/(\gamma\beta_2)}{\bf 1}_{\{u<1\}}+u^{d_2/(\gamma \beta_2)}{\bf 1}_{\{u\geq 1\}})e^{-u}\,{\rm d}u.$$
Since $t>t/d(x,y)^{\gamma}\geq t^{1-1/\beta_2}$ by assumption,
we obtain \eqref{eq:on-diag}. 
As a consequence of the argument above, we obtain the upper bound of \eqref{eq:heat-sub}. 

We next discuss the lower bound of \eqref{eq:heat-sub}. 
Suppose first that $0<t<1$ and $d(x,y)<1$. 
Since $t^{1/\gamma}\vee d(x,y)^{\beta_1}\leq 1\vee d(x,y)$, 
we have
\begin{equation}\label{eq:change}
\begin{split}
&\int_0^{t^{1/\gamma}\vee d(x,y)^{\beta_1}}p(s,x,y)\pi_t(s)\,{\rm d}s\\
&\gtrsim \int_0^{t^{1/\gamma}\vee d(x,y)^{\beta_1}}\frac{1}{V(x,s^{1/\beta_1})}
\exp\left\{-\left(C_1\frac{d(x,y)^{\beta_1}}{s}\right)^{1/(\beta_1-1)}\right\}\pi_t(s)\,{\rm d}s\\
&\geq \frac{1}{V(x, t^{1/\gamma\beta_1}\vee d(x,y))}
\int_0^{t^{1/\gamma}\vee d(x,y)^{\beta_1}} 
\exp\left\{-\left(C_1\frac{d(x,y)^{\beta_1}}{s}\right)^{1/(\beta_1-1)}\right\}
\pi_t(s)\,{\rm d}s.
\end{split}
\end{equation}
If $d(x,y)^{\beta_1}\leq 2t^{1/\gamma}$, 
then by \eqref{eq:volume-growth}, \eqref{eq:scale} and the change of variables $s=t^{1/\gamma}u$, 
the last expression of (\ref{eq:change}) is greater than 
\begin{equation*}
\begin{split}
&\frac{1}{t^{1/\gamma}}\frac{C}{V(x, t^{1/(\gamma\beta_1)})}
\int_0^{t^{1/\gamma}} 
\exp\left\{-\left(C_1\frac{d(x,y)^{\beta_1}}{s}\right)^{1/(\beta_1-1)}\right\}
\pi_1\left(\frac{s}{t^{1/\gamma}}\right)\,{\rm d}s\\
&=\frac{C}{V(x, t^{1/(\gamma\beta_1)})}
\int_0^1 
\exp\left\{-\left(C_1\frac{d(x,y)^{\beta_1}}{t^{1/\gamma}u}\right)^{1/(\beta_1-1)}\right\}
\pi_1(u)\,{\rm d}u\\
&\geq \frac{C}{V(x, t^{1/(\gamma\beta_1)})}
\int_0^1 
\exp\left\{-\left(\frac{2C_1}{u}\right)^{1/(\beta_1-1)}\right\}
\pi_1(u)\,{\rm d}u\asymp \frac{1}{V(x, t^{1/(\gamma\beta_1)})}.
\end{split}
\end{equation*}
On the other hand, if $d(x,y)^{\beta_1}>2t^{1/\gamma}$, 
then we see by \eqref{eq:lower-pi} that 
the last expression of (\ref{eq:change}) is greater than 
\begin{equation*}
\begin{split}
&\frac{1}{V(x,d(x,y))}
\int_{t^{1/\gamma}}^{d(x,y)^{\beta_1}} 
\exp\left\{-\left(C_1\frac{d(x,y)^{\beta_1}}{s}\right)^{1/(\beta_1-1)}\right\}
\pi_t(s)\,{\rm d}s\\
&\geq \frac{c}{V(x,d(x,y))}
\int_{d(x,y)^{\beta_1}/2}^{d(x,y)^{\beta_1}} 
\exp\left\{-\left(C_1\frac{d(x,y)^{\beta_1}}{s}\right)^{1/(\beta_1-1)}\right\}
\frac{t}{s^{1+\gamma}}\,{\rm d}s\\
&\geq c\exp\left\{-\left(2C_1\right)^{1/(\beta_1-1)}\right\}\frac{t}{V(x,d(x,y))}
\int_{d(x,y)^{\beta_1}/2}^{d(x,y)^{\beta_1}} 
\frac{1}{s^{1+\gamma}}\,{\rm d}s
\asymp \frac{t}{V(x,d(x,y))d(x,y)^{\gamma\beta_1}}.
\end{split}
\end{equation*}
Hence we get the lower bound of \eqref{eq:heat-sub}.

Suppose next that $0<t<1$ and $d(x,y)\geq 1$, or $t\geq 1$. 
Since $t^{1/\gamma}\vee d(x,y)^{\beta_2}\geq 1\vee d(x,y)$, 
it follows by (\ref{eq:lower-pi}) that 
\begin{equation*}
\begin{split}
&\int_{t^{1/\gamma}\vee d(x,y)^{\beta_2}}^{\infty}p(s,x,y)\pi_t(s)\,{\rm d}s\\
&\gtrsim t\int_{t^{1/\gamma}\vee d(x,y)^{\beta_2}}^{\infty}\frac{1}{V(x,s^{1/\beta_2})}
\exp\left\{-\left(C_3\frac{d(x,y)^{\beta_2}}{s}\right)^{1/(\beta_2-1)}\right\}\frac{{\rm d}s}{s^{1+\gamma}}\\
&\gtrsim t\int_{t^{1/\gamma}\vee d(x,y)^{\beta_2}}^{\infty}\frac{1}{V(x,s^{1/\beta_2})}
\frac{{\rm d}s}{s^{1+\gamma}}
\asymp \min\left\{\frac{1}{V(x,t^{1/(\gamma\beta_2)})},\frac{t}{V(x,d(x,y))d(x,y)^{\gamma\beta_2}}\right\}.
\end{split}
\end{equation*}
At the last relation above, we used the fact that for any $p>0$,
$$\int_t^{\infty}\frac{{\rm d}s}{V(x,s^{1/\beta_2})s^{1+p}}\asymp \frac{1}{V(x,t^{1/\beta_2})t^p} \quad (t>0),$$ 
which follows by the same way as in \eqref{eq:int-sum}. 
We thus get the lower bound of (\ref{eq:heat-sub}). 
\bigskip

\noindent
{\bf Acknowledgements} \ 
The author would like to thank Professor Naotaka Kajino 
for informing him the references \cite{BBK06, K13} 
and for drawing his attention to the subordinated diffusion processes 
discussed in  Example \ref{ex:sub-diff}. 
He is indebted to Professor Masayoshi Takeda and Professor Jian Wang 
for their valuable comments on the draft of this paper. 
Thank are also due to the two referees for their valuable suggestions 
on the presentation of this paper.

\address{
Yuichi Shiozawa\\
Graduate School of Natural Science and Technology \\
Department of Environmental and Mathematical Sciences \\
Okayama University \\
Okayama 700-8530,
Japan
}
{shiozawa@ems.okayama-u.ac.jp}
\end{document}